%% file: Algebraic.tex
\documentclass[12pt]{article}
\usepackage{amsmath,amssymb,amsfonts,amscd,theorem}
\input{cmds}

\usepackage[latin1]{inputenc}

\title{Algebraic representations of von Neumann algebras}
\author{\bt{c}
{\sc C. Pierre\/}\\
\noalign{\vskip 6pt}
Institut de Mathématique pure et appliquée\\[-6pt]
Université de Louvain\\[-6pt]
Chemin du Cyclotron, 2\\[-6pt]
B-1348 Louvain-la-Neuve,  Belgium\\
pierre@math.ucl.ac.be\te}
\date{\thispagestyle{empty}}

\newcommand{\NN}{\mathbb{N}\,}
\newcommand{\ZZ}{\mathbb{Z}\,}

\newcommand{\CC}{\mathbb{C}\,}
\newcommand{\RR}{\mathbb{R}\,}
\newcommand{\FF}{\mathbb{F}\,}

\newcommand{\MM}{\mathbb{M}\,}

\def\RL{_{R\times L}}

\def\bpr{\paragraph{Proof:}\quad}

\def\backfl#1#2{\smash{\mathop{\hbox to
12mm{\leftarrowfill}} \limits^{\scriptstyle
#1}_{\scriptstyle #2}}}

\def\To{\begin{CD} @>>>\end{CD}}

\def\End{\operatorname{End}}
\def\Tr{\operatorname{Trace}}
\def\det{\operatorname{det}}

\def\Pra{\operatorname{Pra}}
\def\Im{\operatorname{Im}}
\def\Gal{\operatorname{Gal}}
\def\Rep{\operatorname{Rep}}

\def\Repsp{\operatorname{Repsp}}
\def\FRepsp{\operatorname{FRepsp}}
\def\AdFRepsp{\operatorname{AdFRepsp}}
\def\ELLIP{\operatorname{ELLIP}}

\def\Aut{\operatorname{Aut}}

\def\lr{left (resp. right) }
\def\rl{right (resp. left) }

\def\bbf{\boldmath\bf}

{\theorembodyfont{\rmfamily}
\newtheorem{defi}{Definition}[section]
\newtheorem{defis}[defi]{Definitions}}
\newtheorem{lm}[defi]{Lemma}
\newtheorem{propo}[defi]{Proposition}
\newtheorem{coro}[defi]{Corollary}
{\theorembodyfont{\rmfamily}

}

\newcommand\sub[1]{\stepcounter{defi}\noindent{\bf{\thedefi. #1}}:\;\;}

\begin{document} 

\newenvironment{pr}{{\noindent \em Proof:\/}}{}
\pagestyle{empty}

\null\vfill

\begin{center}
{\LARGE Algebraic representations of von Neumann algebras

}
\vfill

{\sc C. Pierre\/}
\vskip 11pt

Institut de Mathématique pure et appliquée\\
Université de Louvain\\
Chemin du Cyclotron, 2\\
B-1348 Louvain-la-Neuve,  Belgium\\
pierre@math.ucl.ac.be

\vfill
%\begin{flushright} first revised version: April 2000\\
%second revised version: May 2002\\
%first received version: September 26, 1997\end{flushright}
\eject

\null\vfill

{\LARGE Algebraic representations of von Neumann algebras

}
\vfill

{\sc C. Pierre\/}
\end{center}
\vskip 11pt

\vfill

\begin{abstract} An (algebraic) extended bilinear Hilbert semispace $H_a^{\mp}$ is proposed as being the natural
representation space for the algebras of von Neumann.  
This bilinear Hilbert semispace has a well defined structure given by the representation space\linebreak $\Repsp(GL_n(L_{\o v}\times L_v))$ 
of an (algebraic) complete bilinear semigroup $GL_n(L_{\o v}\times L_v)$
over the product of sets of completions characterized by increasing ranks.

This representation space is a $GL_n(L^{(nr)}_{\o v}\times L^{(nr)}_v)$-bisemimodule $M^{(nr)}_R\otimes M^{(nr)}_L$~, decomposing into subbisemimodules according to the pseudounramified or pseudoramified conjugacy classes of $GL_n(L^{(nr)}_{\o v}\times L^{(nr)}_v)$~, and is in one-to-one correspondence with its cuspidal representation according to the Langlands global program.

In this context,  towers of von Neumann  subbisemialgebras on graded bilinear Hilbert subsemispaces, of which structures are these subbisemimodules, are constructed algebraically which allows to envisage the classification of the factors of von Neumann from an algebraic point of view.
\end{abstract} \vfill\vfill\eject

\pagestyle{myheadings}
\setcounter{page}{1}

\section*{Introduction}

The first essential step of this paper consists in building up a bilinear
mathematical frame for the representations of the von Neumann algebras in
such a way that the most convenient representation space be essentially an
extended bilinear Hilbert semispace characterized by a non-orthogonal basis.

\vskip 11pt

Considering that the representation space of a von Neumann algebra must be
the enveloping algebra \cite{13} of the Hilbert (semi)module on which this von
Neumann algebra acts, an extended bilinear Hilbert semispace is then proposed whose
Hilbert bisemimodule constitutes the searched enveloping semialgebra \cite{30}: this
constitutes the content of chapter 1 \cite{29}.
\vskip 11pt

In this perspective, an algebraic (real) extended bilinear Hilbert
semispace $H^{\pm}_a$ and an analytic (complex) extended bilinear
Hilbert semispace $H^{\pm}_{h}$ are constructed and proved to be the
natural representation spaces for the algebras of elliptic operators. In
this context, semialgebras and bisemialgebras of von Neumann on the spaces
$H^{\pm}_a$ and $H^{\pm}_{h}$ are introduced
according to the general
treatment of semistructures and  bisemistructures introduced in \cite{30}.
\vskip 11pt

The generation of algebraic bilinear Hilbert semispaces is related to the bilinear Eisenstein cohomology 
which constitutes the  algebraic pillar of the bilinear global program of Langlands 
introduced in \cite{29}.  More concretely, we are interested in the representation space
$\Repsp(GL_n (L_{\o v}\times L_v ))$ of a bilinear general semigroup over the product 
$(L_{\o v}\times L_v )$ of sets of pseudoramified real completions, at infinite archimedean places,  whose ranks (or degrees) are integers modulo $N$ in such a way that:
\Bi
\item $GL_n(L_{\o v}\times L_v )$ has the Gauss bilinear  decomposition into the product of subgroups 
of diagonal matrices by the subgroups of upper and lower unitriangular matrices;

\item $GL_n(L_{\o v}\times L_v )=T_n^t(L_{\o v})\times T_n(L_v )$ has for representation space the tensor 
product $(M_R\otimes M_L)$ of a right $T_n^t(L_{\o v})$-semimodule $M_R$ by a left 
$T_n(L_v )$-semimodule $M_L$ 
such that $M_L$ (resp. $M_R$~) decomposes into $T_n(L_{v_i})$-subsemimodules $M_{v_i}$ 
(resp. $T_n^t(L_{\o v_i})$-subsemimodules $M_{\o v_i}$~) according to  the \lr archimedean places $v_i$ 
(resp. $\o v_i$~) of $L_v $ (resp. $L_{\o v}$~) and so that the set of \lr subsemimodules 	 
$M_{v_i}$ (resp. $M_{\o v_i}$~) corresponds to the set of \lr conjugacy classes of 
$GL_n(L_{\o v}\times L_v )$~.
\Ei
\vskip 11pt

The  bilinear Eisenstein cohomology (semi)group is the cohomology of the Shimura bisemivariety 
given by
\[\partial \o S_{G\RL}= P_n(L_{\o v^1}\times L_{v^1}) \setminus GL_n(L^+_{R}\times L^+_L )
\big/GL_n((\ZZ\big/ N\ \ZZ)^2)\]
where
\Bi
\item $P_n({L_{\o v^1}}\times {L_{v^1}})  $ is a bilinear parabolic subsemigroup over the product, right by left, of sets of irreducible real completions having a rank $N$~;
\item $GL_n((\ZZ\big/ N\ \ZZ)^2)$ is a bilinear arithmetic subsemigroup constituting the representation of the tensor product of Hecke operators and having a representation in a Hecke bilattice; 
\item $L^+_R$ and $L^+_L$ are symmetric (real) algebraic (semi)fields.
\Ei
\vskip 11pt

Then, the bilinear Eisenstein cohomology
\[ H^{2j}( \partial \o S_{G\RL} , \widetilde M^{2j}_R\otimes \widetilde M^{2j}_L))\simeq  \Repsp( GL_{2j}(L_{\o v}\times L_v ))\;, \quad 2j\le n\;, \]
of the Shimura bisemivariety $\partial \o S_{G\RL}$ has coefficients in the 
(bisemi)sheaf $\widetilde M^{2j}_R\otimes \widetilde M^{2j}_L$ over the $GL_{2j}(L_{\o v}\times 
L_v )$-bisemimodule $(M^{2j}_R\otimes M^{2j}_L)$ and is in bijection with the representation space of the complete bilinear algebraic semigroup $GL_{2j}(L_{\o v}\times L_v )$~.

Furthermore, the complete reducibility of $\Repsp(GL_{2n}(L_{\o v}\times L_v ))$ induces the decomposition of the   bilinear Eisenstein cohomology into (irreducible) two-dimensional bilinear Eisenstein cohomologies.
\vskip 11pt

On the other hand, the analytic pillar of the global program of Langlands is given by the 
cuspidal representation of the coefficients of the bilinear Eisenstein cohomology 
in terms of products, right by left, of global elliptic semimodules which are (truncated) Fourier 
series over $\rit$ whose number of terms corresponds to the number of conjugacy classes of the general bilinear
 semigroup $GL_{2j}(L_{\o v}\times L_v )$~.

The Eisenstein and analytic de Rham cohomologies are considered and recalled to be isomorphic 
\cite{21} from which it results that (bi)semialgebras of von Neumann on the algebraic and analytic bilinear 
Hilbert semispaces $H^{\pm}_a$ and $H^{\pm}_h$ are isomorphic:
\[ \MM_{{R,L}\atop\RL}(H^{\pm}_a) \simeq \MM_{{R,L}\atop\RL}(H^{\mp}_h) \;.\]
The action of a (differential) bioperator $(T^D_R\otimes T^D_L)\in \MM\RL(H^{\mp}_a) $ of rank $(m\times m)$~,\linebreak with $m\le n$~, (associated with a principal $GL_m(\rit\times \rit)$-bundle)
 on the $(n\times n)$-dimensional  (bisemi)sheaf $(\widetilde M_R\otimes  \widetilde M_L)$ consists in mapping $(\widetilde M_R\otimes \widetilde M_L)$ into the corresponding   (bisemi)\-sheaf $(\widetilde M_{R_{n[m]}}\otimes \widetilde M _{L_{n[m]}})$ shifted into $(m\times m)$ dimensions such that $(\widetilde M _{R_{n[m]}}\otimes \widetilde M_{L_{n[m]}})$ decomposes into   subbisemisheaves according to:
\Bi
\item the pseudoramified conjugacy classes $g_R(i)\times g_L(i)$~, $1\le i\le q$~, of $GL_n(L_{\o v}\times L_v)$  where ``~$i$~'' denotes a global residue degree;
\item[or \textbullet] according to the pseudounramified conjugacy classes $\gamma_R(i)\times \gamma_L(i)$~, $1\le i\le q$~, of the pseudounramifed bilinear semigroup $GL_n(L^{nr}_{\o v}\times L^{nr}_v)$  over sets $L^{nr}_v$ and $L^{nr}_{\o v}$ of pseudounramified completions;
\Ei
 in such a way that:
\Bi
\item $(\widetilde M_{R_{n[m]}}\otimes \widetilde M_{L_{n[m]}})$ be the coefficient system of the shifted bilinear 
Eisenstein cohomology $H^{2j-2k}(\partial\o S_{G_{R\times L;n[m]}},\widetilde M^{2j}_{R_{2j[2k]}}\otimes 
\widetilde M^{2j}_{L_{2j[2k]}})$ where
\begin{align*}
\partial \o S_{G_{R\times L;n[m]}} =&P_{n[m]}((L_{\o v ^1}\otimes \rit)\times (L_{v^1}\otimes 
\rit))\\
&\quad \setminus GL_{n[m]}((L^+_{R}\otimes \rit)\times (L^+_L \otimes \rit))\Big/GL_{n[m]}(\ZZ\big/N\ \ZZ)^2\otimes \rit^2)\end{align*}
is the shifted Shimura bisemivariety;
\item $(\widetilde M^{2j}_{R_{2j[2k]}}\otimes \widetilde M^{2j}_{L_{2j[2k]}})$ decomposes into shifted  subbisemisheaves according to the pseudoramified or pseudounramified conjugacy classes of $GL_{2j}(L^{(nr)}_{\o v}\times L^{(nr)}_v)$ in such a way that the pseudoramified conjugacy classes correspond to the cosets of $GL_{n[m]}((L^+_{R}\otimes \rit)\times (L^+_L \otimes \rit))\big/ GL_{n[m]}(\ZZ\big/N\ \ZZ)^2\otimes\rit^2)$~.
\Ei

As in the unshifted case, the shifted bilinear Eisenstein cohomology decomposes into direct sum of completely irreducible orthogonal or nonorthogonal shifted bilinear Eisenstein cohomologies.
\vskip 11pt

Taking into account the decomposition of the complete algebraic (resp. analytic) $GL_n(L^{(nr)}_{\o v}\times 
L^{(nr)}_v )$-bisemimodule $(M^{(nr)}_R\otimes M^{(nr)}_L)$ (resp. $(M^{s,(nr)}_R\otimes M^{s,(nr)}_L)$~) into subbisemimodules 
according to its pseudounramified or pseudoramified conjugacy classes, the complete algebraic 
(resp. analytic) extended bilinear Hilbert semispace $H^{\mp,(nr)}_ a$ (resp. $H^{\mp,(nr)}_h$~) also decomposes 
   into bilinear subsemispaces $H^{\mp,(nr)}_a(i)$ 
(resp. $H^{\mp,(nr)}_h(i)$~), $1\le i\le q$, or according to sums of bilinear subsemispaces:
\begin{alignat*}{3}
H^{\mp,nr}_a\{i\} &= \txt\bigoplus\limits^i_{\nu=1}H^{\mp,nr}_a(\nu) \qquad &\text{(resp.} \quad 
H^{\mp,nr}_h\{i\} =& \txt\bigoplus\limits^i_{\nu=1}H^{\mp,nr}_h(\nu)\ ), \tag*{$1\le i\le q$~,}\\[11pt]
\text{or} \quad
H^{\mp}_a\{i\} &= \txt\bigoplus\limits^i_{j=1}H^{\mp}_a(j) \qquad &\text{(resp.} \quad 
H^{\mp}_h\{i\} =& \txt\bigoplus\limits^i_{j=1}H^{\mp}_h(j)\ ).\end{alignat*}
So, towers of sums of embedded bilinear Hilbert subsemispaces
\begin{align*}
H^{\mp,nr}_a\{1\} \subset \cdots \subset & H^{\mp,nr}_a\{i\} \subset \cdots \subset H^{\mp,nr}_a\{q\}\;, \\ 
H^{\mp}_a\{1\} \subset \cdots \subset & H^{\mp}_a\{i\} \subset \cdots \subset H^{\mp}_a\{q\}\;, \\[11pt] 
\text{(resp.} \quad
H^{\mp,nr}_h\{1\} \subset \cdots \subset & H^{\mp,nr}_h\{i\} \subset \cdots \subset H^{\mp,nr}_h\{q\}\;, \\ 
H^{\mp}_h\{1\} \subset \cdots \subset & H^{\mp}_h\{i\} \subset \cdots \subset H^{\mp}_h\{q\}\;), \end{align*}
can be constructed, leading to consider that these bilinear Hilbert semispaces are ``solvable'' and graded.
\vskip 11pt

And towers of sums of pseudounramified or pseudoramified von Neumann sub(bi)\-semialgebras  can be generated according to:
\[ \MM
_{\left\{ {R,L}  \atop \RL\right.}
(H^{\mp,(nr)}_a\{1\})\subset \cdots \subset 
\MM
_{\left\{ {R,L}  \atop \RL\right.}
(H^{\mp,(nr)}_a\{i\})\subset \cdots \subset
\MM
_{\left\{ {R,L}  \atop \RL\right.}
\left(H^{\mp,(nr)}_a \{ q\}\right)\;.\]
Then, the discrete spectrum $\sigma (T^D_R\otimes T^D_L)$ of a (differential) bioperator 
$(T^D_R\otimes T^D_L)\in\linebreak \MM\RL(H^{\mp,(nr)}_a)$ is obtained throughout the morphism from the von Neumann bisemialgebra $\MM\RL(H^{\mp,(nr)}_a)$ to the set of von Neumann subbisemialgebras $[\MM\RL(\Hs^{\mp,(nr)}_a\{i\}]_i$ defined on the set of pseudounramified or pseudoramified bilinear Hilbert subsemispaces $\Hs^{\mp,(nr)}_a\{i\}$ characterized by a diagonal metric associated with an orthonormal bilinear basis.

If the cuspidal representation space of the $GL_n(L^{(nr)}_{\o v}\times L^{(nr)}_v )$-bisemimodule $(M^{(nr)}_R\otimes M^{(nr)}_L)$ is taken into account, the corresponding set of eigenbifunctions of the differential bioperator $(T^D_R\otimes T^D_L)$  is given, according to the Langlands program, by the global elliptic subbisemimodules which are products, right by left, of (truncated) Fourier series (over $\rit$) whose number of terms correspond to the number of archimedean places associated with the considered algebraic intermediate finite number (semi)fields.
\vskip 11pt

In this context, the classification of the factors of von Neumann can be envisaged from the algebraic frame developed in this paper.

In correspondence with the introduction of {\bf pseudoramified\/}  bilinear Hilbert semispaces $H^{\pm}_a$ and of towers of embedded bilinear Hilbert subsemispaces, {\bf pseudounramified\/} bilinear Hilbert semispaces $H^{\rm nr}_a$ can be defined as well as towers of embedded bilinear pseudounramified Hilbert subsemispaces.

So, if ``$i$~'' labels an 
algebraic intermediate (semi)field or the associated archimedean completion, 
$\MM _{R,L}
(H^{\rm nr}_a(i))$ will refer to a factor of type ${\rm I}_i$ while if ``~$j$~'' denotes an algebraic internal dimension, 
$\MM _{R,L}
(H^{\mp,{\rm in}}_a(j))$~, $1\le j\le N$~,  will be a 
hyperfinite subfactor of type ${\rm I}{\rm I}_{1_j}$ \cite{23}, \cite{24}, where $N$ is the order of a global inertia subgroup.

So, our main proposition can finally be stated as follows~:\Be
\item On the pseudounramified bilinear Hilbert semispace $H^{\rm nr}_a$~, there are $q$ factors of type ${\rm I}_i$~, 
$1\le i\le q\le \infty$ where ``~$i$~'' denotes a global residue degree.

\item On the  bilinear Hilbert subsemispace $H^{\rm in}_a[L_{\o v_1}\times L_{v_1}]$ 
restricted to the representation of the bilinear parabolic subsemigroup $P_n(L_{\o v^1}\times 
L_{v^1})$~, there are $N$ subfactors of type ${\rm I}{\rm I}_{1_j}$~, where $j$ denotes an internal 
algebraic dimension.

The upper subfactor ${\rm I}{\rm I}_{1_N}$ is the hyperfinite factor ${\rm I}{\rm I}_1$~.

\item On the tensor products $H^{\text{nr}}_a(i)\otimes H^{\rm in}_a(N)$~, there are $q$ pseudoramified factors of type ${\rm II}_\infty$~, $1\le i\le q\le \infty$~, noted $\MM_{R,L}(H^{\text{nr}}_a(i)\otimes H^{\rm in}_a(N))$ where $i$ denotes a global residue degree.
\item On the tensor products $H^{\rm nr}_a(\infty )\otimes H^{\mp}_a(i)$~, $1\le i\le N$~, the factors of type  ${\rm I}{\rm I}_\infty $ are defined.
\Ee
\vskip 11pt

\section{Bilinear semigroups and bilinear Hilbert semispaces}

The aim of this chapter is to introduce a sufficiently general
mathematical frame for the representations of the von Neumann algebras. 
As the ``representation" of a $k$-algebra $M_L$ over a  number field $k$ of
characteristic zero proceeds from its enveloping algebra, the most natural
representation space for the von Neumann algebras will be an extended
Hilbert semispace of bilinear nature which must then correspond to the
representation space of the $k$-algebra $M_L$ in a linear Hilbert space
$\Hs$~.

If the representation space of a von Neumann algebra is assumed to be non commutative,
its (alge\-braic)-geometric structure will then be of Riemann type and
composed of the tensor product of a pair of faithfully projective isomorphic $k$-semimodules
leading to an extended bilinear Hilbert semispace by projection of one of
these semimodules on its copy.
\vskip 11pt

\paragraph{Notations:} \bt[t]{lll}
$R,L$ &means& ``~$R$~" or ``~$L$~" for ``right'' or ``left'';\\ $\times_{(D)}$ &means & a diagonal
(~$\times_D$~) or complete (~$\times$~) product.\te
\vskip 11pt

\begin{defi}\quad {\bf Enveloping algebra:}\quad  
Let $M_{R,L}$ be a $k$-algebra considered as a finitely generated,
projective and faithful right (resp. left) $k$-module.  Its enveloping
algebra is given by $M^e=M_R\otimes_kM_L$ where $M_R$ (resp. $M_L$~) is a
right (resp. left) $k$-module viewed as the opposite algebra of $M_L$
(resp. $M_R$~) \cite{13}.  If the homomorphism $E_{h_{R,L}}:M^e\to {\rm
End}_k(M_{R,L})$ is an isomorphism, then the $k$-algebra $M_{R,L}$ is
called an Azumaya algebra.

If $M_{R,L}$ is a faithfully projective right (resp. left) $k$-module of
dimension $n$~, then $M_{R,L}\simeq k^n$
 and we have that \cite{6}, \cite{18}, \cite{35}:
\[ M^e\simeq {\rm End}_k(M_{R,L})\simeq {\rm End}_k(k^n)\simeq M_n(k)\]
where $M_n(k)$ is the ring of matrices of order $n$ over $k$~.

The homomorphism $E_{R,L}:M_{R,L}\to M_n(k)$ is called a $n$-dimensional
representation of $M_{R,L}$ \cite{2}.
\end{defi}
\vskip 11pt

\begin{defi}\quad {\bf Symmetric algebraic extension field:}\quad   
Let $k$ be  a number field of characteristic $0$ and $L^+$ (resp. $L$~) denote a finite real (resp. complex) extension of $k$~.  A real (resp. complex) algebraic
extension field $L^+$ (resp. $L$~) will be said symmetric if it is composed of the set of
positive  (resp. complex) simple roots, noted $L^+_L$ (resp. $L_L$~), in one-to-one correspondence with the
set of negative (resp. complex conjugate) simple roots, noted $L^+_R$ (resp. $L_R$~), such that to each positive (resp. complex)
simple root $x^+_L\in L_L^+$ (resp. $x_L\in L_L$~) corresponds a symmetric negative (resp. complex conjugate) simple root
$x^+_R\in L^+_R$ (resp. $x_R\in L_R$~). Geometrically, $L_L$ is then localized in the upper
halfspace and $L_R$ in the lower half space.  $L^+_L$ (resp. $L_L$~) and $L^+_R$ (resp. $L_R$~) are then
respectively left and right semifields, i.e. commutative division left and
right semirings.

$L_L$ and $L_R$ are semirings because they are abelian semigroups with
respect to the addition and are endowed with associative multiplication
and distributive laws.
\end{defi}
\vskip 11pt

\begin{defi}\quad {\bf Completions associated with finite algebraic extensions:}\quad   
The equivalence classes of the real completions of $L^+_L$ (resp. $L^+_R$~), obtained by an isomorphism of compactification of the corresponding extensions, are the \lr infinite places of $L^+_L$ (resp. $L^+_R$~) and are noted $v=\{v_1,\cdots,v_i,\cdots,v_q\}$ (resp. $\o v=\{\o v_1,\cdots,\o v_i,\cdots,\o v_q\}$~).

Similarly, the equivalence classes of the complex completions of $L_L$ (resp. $L_R$~), obtained by an isomorphism of compactification of corresponding finite extensions, are the \lr infinite complex places of $L_L$ (resp. $L_R$~) and are noted 
$\omega =\{\omega _1,\dots,\omega _i,\dots,\omega _q\}$
(resp. $\o\omega =\{\o\omega _1,\dots,\o\omega _i,\dots,\o\omega _q\}$~).

Let $L_{v_i}$ (resp. $L_{\o v_i}$~) denote $i$-th basic real completion corresponding to the $i$-th \lr pseudoramified algebraic extension $L^+_{L_i}$ (resp. $L^+_{R_i}$~) of $k$ and associated to the \lr place $v_i$ (resp. $\o v_i$~).  The other equivalent completions of $v_i$ (resp. $\o v_i$~) are noted $L{v_i,m_i}$ (resp. $L_{\o v_i,m_i}$~), where $m_i\in\NN$~, $m_i>0$~, are increasing integers.\\
(~$m_i=0$ refers to the basic completion $L_{v_i}$ (resp. $L_{\o v_i}$~) ).

It is assumed that the \lr pseudoramified completions $L_{v_i,m_i}$ (resp. $L_{\o v_i,m_i}$~) are generated from an irreducible completion $L_{v_i^1}$ (resp. $L_{\o v_i^1}$~) having a rank or degree equal to $N$~.

Then, the rank of the pseudoramified completions $L_{v_i}$ (resp. $L_{\o v_i}$~) and $L_{v_i,m_i}$ (resp. $L_{\o v_i,m_i}$~), corresponding to the degree of extension of the associated extension, is given by an integer modulo $N$ according to:
\begin{align*}
n_{i_L} &= [L_{v_i,m_i}:k]=*+f_{v_i}\cdot N\simeq i\cdot N\\
\text{(resp.}\quad n_{i_R} &= [L_{\o v_i,m_i}:k]=*+f_{\o v_i}\cdot N\simeq i\cdot N\ )\end{align*}
where
\Bi
\item $*$ denotes an integer inferior to $N$~;
\item $f_{v_i}$ (resp. $f_{\o v_i}$~), called a global class residue degree, is the degree of the corresponding pseudounramified completions $L^{nr}_{v_i,m_i}$ (resp. $L^{nr}_{\o v_i,m_i}$~) given by
\[[ L^{nr}_{v_i,m_i}:k]=f_{v_i}=i \qquad 
\text{(resp.}\quad L^{nr}_{\o v_i,m_i}:k]=f_{\o v_i}=i\ )\]
So, the ranks or degrees of the pseudoramified completions $L_{v_i,m_i}$ (resp. $L_{\o v_i,m_i}$~), $1\le i \le q$~, 
are integers modulo $N$~, $\ZZ\big/N\ \ZZ$~.
\Ei

Remark that the integer $\sup(m_i)$ is interpreted as the multiplicity of the place $v_i$ (resp. $\o v_i$~).

As the rank $n_{i_L}$ (resp. $n_{i_R}$~) of the completion $L_{v_i,m_i}$ (resp. $L_{\o v_i,m_i}$~) is assumed to be a multiple of the integer $N$~, which is the rank of the irreducible subcompletion $L_{v^1_i}$ (resp. $L_{\o v^1_i}$~), the completion $L_{v_i,m_i}$ (resp. $L_{\o v_i,m_i}$~)  will be cut into a set of $i$ equivalent subcompletions
$L_{v^{i'}_i,m_i}$ (resp. $L_{\o v^{i'}_i,m_i}$~), $1\le i'\le i$~, of rank $N$~.

Finally, let
\begin{align*}
L_v &= \{L_{v_1},\cdots,L_{v_{i,m_i}},\cdots,L_{v_{q,m_q}}\}\\
\text{(resp.} \quad
L_{\o v}&= \{L_{\o v_1},\cdots,L_{\o v_{i,m_i}},\cdots,L_{\o v_{q,m_q}}\}\ )\end{align*}
denote the set of real pseudoramified completions of $L^+_L$ (resp. $L^+_R$~) with
\[ L_{v_\oplus}=\txt\bigoplus\limits_i\bigoplus\limits_{m_i} L_{v_{i,m_i}}
\qquad \text{(resp.} \quad
L_{\o v_\oplus}=\txt\bigoplus\limits_i\bigoplus\limits_{m_i} L_{\o v_{i,m_i}}\ )\]
be their direct sum and let
\begin{align*}
L^{nr}_v &= \{L^{nr}_{v_1},\cdots,L^{nr}_{v_{i,m_i}},\cdots,L^{nr}_{v_{q,m_q}}\}\\
\text{(resp.} \quad
L^{nr}_{\o v}&= \{L^{nr}_{\o v_1},\cdots,L^{nr}_{\o v_{i,m_i}},\cdots,L^{nr}_{\o v_{q,m_q}}\}\ )\end{align*}
denote the corresponding set of real pseudounramified completions.

Similarly, let
\begin{align*}
L_\omega  &= \{L_{\omega _1},\cdots,L_{\omega _{i,m_i}},\cdots,L_{\omega _{q,m_q}}\}\\
\text{(resp.} \quad
L_{\o \omega }&= \{L_{\o \omega _1},\cdots,L_{\o \omega _{i,m_i}},\cdots,L_{\o \omega _{q,m_q}}\}\ )\end{align*}
 denote the set of complex pseudoramified completions of $L_L$ (resp. $L_R$~) in such a way that the set $L_v$ (resp. $L_{\o v}$~) of real completions covers the corresponding set $L_\omega $ (resp. $L_{\o\omega }$~) of complex completions \cite{29}.
\end{defi}
\vskip 11pt

\begin{defi}\quad {\bf Galois subgroups and inertia subgroups:}\quad   
Let $\Gal(L^+_{L_i}/k)$ (resp. $\Gal(L^+_{R_i}/k)$~) be the Galois subgroup of the pseudoramified extension $L^+_{L_i}$ (resp. $L^+_{R_i}$~) and let $\Gal(L^{nr,+}_{L_i}/k)$ (resp. $\Gal(L^{nr,+}_{R_i}/k)$~) denote the Galois subgroup of the corresponding pseudounramified extension $L^{nr,+}_{L_i}$
(resp. $L^{nr,+}_{R_i}$~).

If ${\rm I}_{L^+_{L_i}}$ (resp. ${\rm I}_{L^+_{R_i}}$~), denoting the global inertia subgroup of $\Gal(L^+_{L_i}/k)$ 
(resp.\linebreak $\Gal(L^+_{R_i}/k)$~), is the group of Galois automorphisms of the irreducible extension $L^+_{L^1_i}$ (resp. $L^+_{R^1_i}$~) or the group of Galois inner automorphisms,
then we have that
\begin{align*}
\Gal(L^+_{L_i}/k)\big/{\rm I}_{L^+_{L_i}}&=\Gal (L^{nr,+}_{L_i}/k)\\[11pt]
 \text{(resp.} \quad 
\Gal(L^+_{R_i}/k)\big/{\rm I}_{L^+_{R_i}}&=\Gal (L^{nr,+}_{R_i}/k)\ )\end{align*}
such that the exact sequence:
\begin{alignat*}{7}
0 & \To {\rm I}_{L^+_{L_i}} & \To & \Gal (L^+_{L_i}/k) & \To & \Gal (L^{nr,+}_{L_{L_i}}/k) &\To & 1\\
\text{(resp.} \quad 
0 & \To {\rm I}_{L^+_{R_i}} & \To & \Gal (L^+_{R_i}/k) & \To & \Gal (L^{nr,+}_{L_{R_i}}/k) 
&\To & 1\ )\end{alignat*}
has kernel given by the global inertia subgroup ${\rm I}_{L^+_{L_i}}$ (resp. ${\rm I}_{L^+_{R_i}}$~) associated to the place $v_i$ (resp. $\o v_i$~).

If $m_i$ refers to the multiplicity of the left and right places $v_i $ and $\o v_i$~, then the \lr Galois group can be decomposed according to:
\begin{align*}
\Gal(L^+_L/k)&= \txt \bigoplus\limits^q_{i=1} \bigoplus\limits_{m_i} \Gal (L^+_{L_i,m_i}/k)\\[11pt]
 \text{(resp.} \quad 
\Gal(L^+_R/k)&= \txt \bigoplus\limits^q_{i=1} \bigoplus\limits_{m_i} \Gal (L^+_{R_i,m_i}/k) \ ).
\end{align*}
\end{defi}
\vskip 11pt

\sub{Representation of the bilinear general semigroup} \quad Let ${L_v}$ (resp. ${L_{\o v}}$~) be the set of pseudoramified real completions of $L^+_L$ (resp. $L^+_R$~).  Then, a bilinear general (or complete algebraic) semigroup over the product $L_{\o v}\times L_v$  can be defined as the product of the (semi)group $T^t_n(L_{\o v})$ of  lower triangular matrices of order $n$ over $L_{\o v}$ by the (semi)group $T_n(L_{v})$ of upper triangular matrices of order $n$ over $L_v$ according to \cite{29}:
\[GL_n(L_{\o v}\times L_v) = T^t_n(L_{\o v})\times T_n(L_v)\]
such that:
\Bean
\item $GL_n(L_{\o v}\times L_v)$ has the bilinear Gauss decomposition:
\[ GL_n(L_{\o v}\times L_v) =[(D_n(L_{\o v})\times D_n(L_v)][UT_n(L_v)\times UT^t_n(L_{\o v})]\]
where
\Bi
\item $D_n(\cdot)$ is the subgroup of diagonal matrices;
\item $UT_n(\cdot)$ is the subgroup of unitriangular matrices.
\Ei

\item $GL_n(L_{\o v}\times L_v)$ has for  representation space $\Repsp(GL_n(L_{\o v}\times L_v)) $ given by the tensor product $M_R\otimes M_L$ of a right $T^t_n(L_{\o v})$-semimodule $M_R$ localized in the upper half space by a left $T_n(L_v)$-semimodule $M_L$ localized in the lower half space.

\item the \lr conjugacy classes of $T_n(L_v)$ (resp. $T^t_n(L_{\o v})$~) 
correspond to the \lr places $v_i$ (resp. $\o v_i$~) of $L_v$ (resp. $L_{\o v}$~).
\Ee

Similarly, $GL_n(L^{nr}_{\o v}\times L^{nr}_v)$ has for  representation space $\Repsp(GL_n(L^{nr}_{\o v}\times L^{nr}_v)) $ given by the tensor product $M^{nr}_R\otimes M^{nr}_L$ of a right pseudounramified $T^t_n(L^{nr}_{\o v})$-semimodule $M^{nr}_R$  by its left equivalent $T_n(L^{nr}_v)$-semimodule $M^{nr}_L$~.

Considering complete bilinear algebraic (semi)groups is justified by the fact that they ``cover'' their ``linear'' equivalents.  Indeed, it was proved in \cite{29} that a linear complete algebraic group $GL_n(L_{\o v-v})$~, with entries in $L_{\o v-v}\equiv L_{\o v}\cup L_v$ and representation space given by a vectorial space $V$ of dimension $n^2$~, is covered by the bilinear complete algebraic semigroup $GL_n(L_{\o v}\times L_v)$~, having as representation space the $GL_n(L_{\o v}\times L_v)$-bisemimodule $M_R\otimes M_L$~, at the conditions given in \cite{29}.

On the other hand, let $M_{L_\oplus}$ (resp. $M_{R_\oplus}$~) denote the representation space of $T_n(L_{v_\oplus})$ (resp. $T^t_n(L_{\o v_\oplus})$~) with entries in the sum $L_{v_\oplus}$ (resp. $L_{\o v_\oplus}$~) of real pseudoramified completions $L_{v_{i,m_i}}$ (resp. $L_{\o v_{i,m_i}}$~).

Then, $M_{L_\oplus}$ (resp. $M_{R_\oplus}$~) is homomorphic to $M_L$ (resp. $M_R$~) and decomposes into the direct sum of 
$T^t_n(L_{v_i})$-subsemimodules $M_{v_i}$ (resp.
$T^t_n(L_{\o v_i})$-subsemimodules $M_{\o v_i}$~) according to:
\[ M_{L_\oplus} =\txt\bigoplus\limits^q_{i=1} \bigoplus\limits_{m_i} M_{v_i;m_i} \qquad \text{(resp.} \quad 
M_{R_\oplus} =\txt\bigoplus\limits^q_{i=1} \bigoplus\limits_{m_i} M_{\o v_i;m_i}\ )\]
such that:
\Bean
\item each $T_n(L_{v_i})$-subsemimodule $M_{v_i}$ (resp. $T^t_n(L_{\o v_i})$-subsemimodule 
$M_{\o v_i}$~) of dimension $n$ constitutes a representative of the $i$-th conjugacy class of 
$T_n(L_v)$ (resp. $T^t_n(L_{\o v})$~).

\item $M_{v_i}$ and $M_{\o v_i}$~, $1\le i\le q$~, has a rank given by:
\[ n_i=i^n\cdot N^n=f^n_{v_i}\cdot N^n\;.\]
\Ee
%\end{rien}
\vskip 11pt

\begin{defi}\quad {\bf Bisemimodules:}\quad  
The bilinear tensor product between the right\linebreak $T^t_n(L_{\o v})$-semimodule $M_R$ and the left $T_n(L_v)$-semimodule
$M_L$ is given by \cite{30}:
\[T_X:\quad \begin{array}[t]{llll}
\{M_R,M_L\}&\To& M_R\otimes M_L\;,&\\
\{x_R,x_L\}&\To & x_R\times x_L\;, & \forall\ x_R\in M_R\ , \;
x_L\in M_L\;,\end{array}\]
so that the pair $\{x_R,x_L\}$ of right and left points be mapped
into the bipoint $x_R\times x_L$ characterized by a Riemanian signature
\cite{18}.  $M_R\otimes M_L$ then is a $GL_n(L_{\o v}\times L_v)$-bisemimodule.

Similarly, the diagonal tensor product between the right and left
semimodules $M_R$ and $M_L$ can be defined by
\[T_{X_D}:\quad \begin{array}[t]{lll}
\{M_R,M_L\}&\To& M_R\otimes_D M_L\;,\\
\{x_R,x_L\}&\To & x_R\times_D x_L\;,\end{array}\]
so that the ``diagonal" bipoint $x_R\times_D x_L$ be characterized by a
diagonal signature which can be Euclidian or not following that the metric
be given by a diagonal unit matrix or by a diagonal matrix having diagonal
elements taking values in the considered field.

$M_R\otimes_D M_L$ then is a $GL_n(L_{\o v}\times_D L_v)$-bisemimodule.
\end{defi}
\vskip 11pt

\begin{defi}\quad {\bf Bisemisheaves of rings:} \quad We want to introduce the set of smooth differentiable (bi)functions on the $GL_n(L_{\o v}\times L_v)$-bisemimodule $M_R\otimes M_L$~, in such a way that these bifunctions are tensor products $\phi _{G_R}(x_{g_R})\otimes \phi _{G_L}(x_{g_L})$ of smooth differentiable right functions $\phi _{G_R}(x_{g_R})$~, $x_{g_R}\in T^t_n(L_{\o v})$~, on $M_R$~, localized in the lower half space by symmetric smooth differentiable left functions 
$\phi _{G_L}(x_{g_L})$~, $x_{g_L}\in T_n(L_v)$~, on $M_L$~, localized in the upper half space.

As $GL_n(L_{\o v}\times L_v)$ is partitioned into conjugacy classes, we have to take into account the bifunctions
$ \phi _{G_{i,m_{i_R}}}(x_{i_R}) \otimes  \phi _{G_{i,m_{i_L}}}(x_{i_L}) $ on the conjugacy class representatives
$ M_{\o v_{i,m_i}}\otimes M_{v_{i,m_i}}$~.  The set of smooth differentiable bifunctions
$\{\phi _{G_{i,m_{i_R}}}(x_{i_R}) \otimes  \phi _{G_{i,m_{i_L}}}(x_{i_L}) \}_{i,m_i}$ on the $GL_n(L_{\o v}\times L_v)$-bisemimodule $M_R\otimes M_L$ is a bisemisheaf of rings noted $\widetilde M_R\otimes \widetilde M_L$ in such a way that this set of differentiable bifunctions are the (bi)sections of $\widetilde M_R\otimes \widetilde M_L$~.

Indeed, $\widetilde M_R$ (resp. $\widetilde M_L$~), having as sections the smooth differentiable functions\linebreak
$\phi _{G_{i,m_{i_R}}}(x_{i_R}) $ (resp. $\phi _{G_{i,m_{i_L}}}(x_{i_L}) $~), is a semisheaf of rings because it is a sheaf of abelian semigroups $\widetilde M_R(x_{i_R})$ (resp. $\widetilde M_L(x_{i_L})$~) for every \rl point $x_{i_R}$ (resp. $x_{i_L}$~) of the topological semispace $M_R=\Repsp(T^t_n(L_{\o v}))$ (resp.  $M_L=\Repsp(T_n(L_v))$~) where $\widetilde M_R(x_{i_R})$ (resp. $\widetilde M_L(x_{i_L})$~) has the structure of a semiring.

The introduction of the bilinear Hilbert semispaces in the next section will concern the bisemisheaf of rings $\widetilde M_R\otimes_{(D)}\widetilde M_L$ as well as the $GL_n(L_{\o v}\times L_v)$-bisemimodule $M_R\otimes_{(D)}M_L$~, on which it is defined, but the developments will only bear on $M_R\otimes_{(D)}M_L$ for the simplicity of the notations.
\end{defi}\vskip 11pt

\begin{defis}\quad {\bf a)\ \ External diagonal bilinear Hilbert semispaces $\Hs^a_L$ and
$\Hs^a_R$:}\quad Let $M_R\otimes_D M_L$ be the diagonal $GL_n(L_{\o v}\times_D L_v)$-bisemimodule.  Consider the projective linear mapping
$p_L:M_R\otimes_D M_L\to M_{R(P)/L}$ projecting the  $T^t_n(L_{\o v})$-semimodule
$M_R$ on the $T_n(L_v)$-semimodule $M_L$~.  $M_{R(P)/L}$ is a bisemimodule
representable locally by the bilinear Hilbert scheme ${\rm
Hilb}_{S_{R(P)}/S_L}$ (case $\widetilde M_{R(P)/L}$~) \cite{xx}.

If $M_{R(P)/L}$ is endowed with an external scalar product $\langle
\phi_P,\psi\rangle$ defined from $M_{R(P)}\times_D M_L$ to $\CC$~,
$\forall\ \phi_P\in M_{R(P)}$~, $\forall\ \psi\in M_L$~, this bisemimodule
$M_{R(P)/L}$ will be called a left external bilinear Hilbert semispace, noted
$\Hs^a_L$~.

Similarly, if we consider the projective linear mapping
$p_R:M_R\otimes_DM_L\to M_{L(P)/R}$ projecting the $T_n(L_v)$-semimodule $M_L$
on the $T^t_n(L_{\o v})$-semimodule $M_R$~, we generate the bisemimodule $M_{L(P)/R}$
representable locally by the bilinear Hilbert scheme ${\rm
Hilb}_{S_{L(P)}/S_R}$~.

Endowing $M_{L(P)/R}$ with an external scalar product from
$M_{L(P)}\times_DM_R$ to $\CC$~, we shall get a right external bilinear
Hilbert semispace noted $\Hs^a_R$~.

Notice that $\Hs^a_L$ and $\Hs^a_R$ are characterized by ortho(normal) basis.
\vskip 11pt

\noindent {\bf b)\ \ Internal diagonal bilinear Hilbert semispaces $\Hs^-_a$ and $\Hs^+_a$:} \quad  Let
$B_L:M_{R(P)}\to M_L$ (resp. $B_R:M_{L(P)}\to M_R$~) be a bijective linear
isometric map from $M_{R(P)}$ (resp. $M_{L(P)}$~) to $M_L$ (resp. $M_R$~)
mapping each covariant element of $M_{R(P)}$ (resp. $M_{L(P)}$~) noted
$M_{L_R}$ (resp. $M_{R_L}$~) into a contravariant element of $M_L$ (resp. $M_R$~).

Then, $B_L$ (resp. $B_R$~) transforms the left (resp. right) external
Hilbert semispace $\Hs_L^a$ (resp. $\Hs_R^a$~) into the left (resp. right)
internal bilinear Hilbert semispace $\Hs^+_a$ (resp. $\Hs^-_a$~) in such a way that
\Bean\item the bielements of $\Hs^+_a$ (resp. $\Hs^-_a$~) are bivectors, i.e.
two confounded vectors;
\item each external scalar product of $\Hs_L^a$ (resp. $\Hs_R^a$~) is
transformed into an internal scalar product defined from $M_{L_R}\times_D
M_L$ (resp. $M_{R_L}\times_D M_R$~) to $\CC$~.
\item $\Hs^+_a$ and $\Hs^-_a$ are characterized by ortho(normal) basis.\Ee
\vskip 11pt

\noindent {\bf c)\ \ Extended external bilinear Hilbert semispaces $H_L^a$ and $H_R^a$:} \quad  If
we consider on the non-Euclidian $GL_n(L_{\o v}\times L_v)$-bisemimodule $M_R\otimes
M_L$ the projective linear mapping $p_L:M_R\otimes M_L\to M_{R(P)/_cL}$
(``~$c$~" for complete), (resp. $p_R:M_R\otimes M_L\to M_{L(P)/_cR}$~) of
the right (resp. left) semimodule $M_R$ (resp. $M_L$~) on the left (resp.
right) semimodule $M_L$ (resp. $M_R$~), we get the non-Euclidian
bisemimodule $M_{R(P)/_cL}$ (resp. $M_{L(P)/_cR}$~).

If we endow $M_{R(P)/_cL}$ (resp. $M_{L(P)/_cR}$~) with a complete external
bilinear form defined from $M_{R(P)}\times M_L$ (resp. $M_{L(P)}\times
M_R$~) to $\CC$~, we get a left (resp. right) extended external bilinear
Hilbert semispace noted $H_L^a$ (resp. $H_R^a$~) characterized by a
non-Euclidian geometry and a non-orthogonal basis.
\vskip 11pt

\noindent {\bf d)\ \ Extended internal bilinear Hilbert semispaces $H^+_a$ and
$H^-_a$:} \quad  The left (resp. right) extended external bilinear Hilbert semispace
$H_L^a$ (resp. $H_R^a$~) can be transformed into the left (resp. right)
extended internal bilinear Hilbert semispace $H^+_a$ (resp. $H^-_a$~) by means of
a bijective linear isometric map $B_L$ (resp. $B_R$~) from $M_{R(P)}$
(resp. $M_{L(P)}$~) into $M_L$ (resp. $M_R$~).

The complete external bilinear form of $H_L^a$ (resp. $H_R^a$~) is then
transformed into a complete internal bilinear form of $H^+_a$ (resp. $H^-_a$~).
\end{defis}
\pagebreak%vskip 11pt

\section{Cohomologies and representation spaces of algebras of operators}

We are interested in the cohomology of compact spaces \cite{8}. So, the most
evident algebraic cohomology of compact spaces is the Eisenstein
cohomology which is based upon the Borel-Serre compactification of the
lattice space attached to an arithmetic group $\Gamma $~.  The
Eisenstein cohomology classes were assumed to be represented by differential forms
which are Eisenstein series \cite{21}, \cite{22}, \cite{32}.
\vskip 11pt

\begin{defi} {\bf The Shimura bisemivariety:\/} Referring to the linear lattice space $X=GL_n(\rit)\big/ GL_n(\ZZ)$~, \cite{6}, \cite{7}, a bilinear complex lattice bisemispace can be introduced by:
\[X_{S\RL}=GL_n(L_R^{(\CC)}\times L_L^{(\CC)})\big/ GL_n((\ZZ/N\ \ZZ)^2)\]
where
\Bi
\item $GL_n((\ZZ/N\ \ZZ)^2)$ is a bilinear arithmetic semigroup over squares of integers modulo $N$~;
\item $GL_n(L_R^{(\CC)}\times L_L^{(\CC)})$ is a bilinear algebraic semigroup with entries in the product $(L_R^{(\CC)}\times L_L^{(\CC)})$ of complex symmetric (semi)fields associated with $(L_R\times L_L)$~.
\Ei
The boundary $\partial \o X_{S\RL}$ of the compactified bisemispace $\o X_{S\RL}$ 
corresponds to the boundary of the Borel-Serre compactification and is given by:
\[ \partial \o X_{S\RL}=GL_n(L^+_{R_d}\times L_{L^+_d})\big/ GL_n(\ZZ/N\ \ZZ)^2)\approx GL_n(L_{\o v}\times L_v)\]
where
$L^+_{R_d}$ and $L^+_{L_d}$ are real compact semifields generated from $L^+_R$ and $L^+_L$~.

The double coset decomposition $\partial \o S_{G\RL}$ of the boundary $\partial \o X_{S\RL}$ of the compactified lattice bisemispace corresponds to a Shimura bisemivariety and is given by:
\[ \partial \o S_{G\RL}=P_n(L_{\o v_1}\times L_{v_1})\setminus GL_n(L^+_{R_d}\times L^+_{L_d}) \big/ GL_n((\ZZ/N\ \ZZ)^2)\approx GL_n(L_{\o v}\times L_v)\]
where
\Bi
\item the subgroup $GL_n((\ZZ/N\ \ZZ)^2)$ constitutes the representation of the coset representatives of the tensor product $T_R(n;q)\otimes T_L(n;q)$ of Hecke operators  \cite{29};

\item $P_n({L_{v_1}})$ is the standard parabolic subsemigroup over the set ${L_{v^1}}= \{L_{v_1^1} ,\cdots,L_{v^1_{i,m_i}},\linebreak \cdots,L_{v^1_{q,m_q}}\}$ of irreducible completions $L_{v^1_{i,m_i}}$ having a rank $N$~. 
$P_n({L_{\o v^1}}\times {L_{v^1}})$ is then a bilinear parabolic subsemigroup constituting the smallest connected pseudoramified normal bilinear subsemigroup of $GL_n(L_{\o v}\times L_v)$ and representing the $n$-fold product ${\rm I}_{L_{\o v_i}}\times {\rm I}_{L_{v_i}}$ of global inertia subgroups.
\Ei

The double coset decomposition $\partial \o S_{G\RL}$~, corresponding to a Shimura bisemivariety and restricted to the lower (resp. upper) half space, becomes:
\begin{align*}
\partial \o S_{G_R} &= P_n({L_{\o v^1}})\setminus T^t_n(L^+_{R_d})\big/ T^t_n(\ZZ/N\ \ZZ)\\
\text{(resp.} \quad 
\partial \o S_{G_L} &= P_n({L_{v^1}})\setminus T_n(L^+_{L_d})\big/ T_n(\ZZ/N\ \ZZ)\ ).\end{align*}
\end{defi} \vskip 11pt

\begin{propo}  \quad The (bi)cosets of the bilinear quotient semigroup $GL_n(L^+_{R_d}\times L^+_{L_d})\big/\linebreak GL_n((\ZZ/N\ \ZZ)^2)$ coincide with the  conjugacy classes of the general bilinear semigroup $GL_n(L_{\o v}\times L_v)$ with respect to the smallest connected pseudoramified normal bilinear subsemigroup given by the bilinear parabolic subsemigroup $P_n({L_{\o v^1}}\times {L_{v^1}})$~.
\end{propo} \vskip 11pt

\paragraph{Sketch of the proof:} \quad According to 1.5, the conjugacy classes of 
$GL_n(L_{\o v}\times L_v)$ are in one-to-one correspondence with the (bi)places of 
$L_{\o v}\times L_v$~.  And, on the other hand, the bilinear subsemigroup $GL_n((\ZZ/N\ \ZZ)^2)$ 
is a  representation of the tensor product of Hecke operators such that the $i$-th (bi)coset 
representative of $GL_n((\ZZ/N\ \ZZ)^2)$ corresponds to the biplace $\o v_i\times v_i$ of $L_{\o v}\times L_v$~.
\epr \vskip 11pt

\begin{propo} \quad The  bilinear cohomology (semi)group of the Shimura bisemivariety
\[\partial \o S_{G\RL} =P_n({L_{\o v^1}}\times {L_{v^1}})\setminus GL_n(L^+_{R_d}\times L^+_{L_d})\big/ GL_n((\ZZ/N\ \ZZ)^2)\]
has its coefficient system given by the  bisemisheaf $( \widetilde M^{2j}_R\otimes \widetilde M^{2j}_L)$ and is given by the bilinear Eisenstein cohomology:
\[H^{2j}(\partial \o S_{G\RL},\widetilde M^{2j}_R\otimes \widetilde M^{2j}_L)\simeq \Repsp(GL_{2j}(L_{\o v}\times L_v))\;, \quad 2j\le r\;,\]
which:
\Bi
\item is in bijection with the representation space $\Repsp(GL_{2j}(L_{\o v}\times L_v))$ of the  bilinear general semigroup $GL_{2j}(L_{\o v}\times L_v)$~;
\item decomposes according to the conjugacy classes of $GL_{2j}(L_{\o v}\times L_{v})$.
\Ei
\end{propo} \vskip 11pt

\bpr 
\Be \item It was demonstrated in \cite{29} that the bilinear Eisenstein cohomology $H^{n}(\partial \o S_{G\RL},\linebreak \widetilde M^{2j}_R\otimes \widetilde M^{2j}_L)$ is in bijection with the representation of the bilinear general semigroup $GL_{2j}(L_{\o v}\times L_v)$~: this results from the fact that the Eisenstein bilinear cohomology can be deduced from the Weil bilinear algebra of the Lie bilinear nilpotent semialgebra.
\item As the bicosets of 
$\partial \o S^{(2j)}_{G\RL} =P_{2j}({L_{\o v^1}}\times {L_{v^1}})\setminus GL_{2j}(L^+_{R_d}\times L^+_{L_d})\big/ GL_{2j}((\ZZ/N\ \ZZ)^2)$
 coincide with the conjugacy classes of $GL_{2j}(L_{\o v}\times L_v)$~, we have that the bilinear Eisenstein cohomology decomposes according to:
\[H^{2j}(\partial \o S_{G\RL},\widetilde M^{2j}_{R_\oplus}\otimes \widetilde M^{2j}_{L_\oplus})\simeq \txt\bigoplus\limits^q_{i=1}
\bigoplus\limits_{m_i}(\widetilde M^{2i}_{\o v_i;m_i}\otimes \widetilde M^{2i}_{v_i;m_i})\;.\]
To each $T_{2j}(L_{v_i})$-subsemimodule $M^{2j}_{v_i}$ (resp. $T^t_{2j}(L_{\o v_i})$-subsemimodule $M^{2j}_{\o v_i}$~) is associated a weight $\lambda _{L_i}$ (resp. $\lambda _{R_i}$~) characterizing the $i$-th \lr Hecke sublattice.  Indeed, there exists the surjective morphism:
\[ i_{M_{L,R}}: \quad M^{2j}_{L,R_\oplus}\To \Lambda _{L,R}\]
from the $T_{2j}(L_v)$-semimodule $M^{2j}_L$ (resp. $T^t_{2j}(L_{\o v})$-semimodule $M^{2j}_R$~) into the\linebreak $T_{2j}(\ZZ/N\ \ZZ)$-semimodule $\Lambda _L$ (resp. $T^t_{2j}(\ZZ/N\ \ZZ)$-semimodule $\Lambda _R$~) which is a \lr Hecke lattice decomposing according to the conjugacy classes of $T_{2j}(L_v)$ (resp. $T^t_{2j}(L_{\o v})$~):
\[ \Lambda _L =\txt\bigoplus\limits^q_{i=1}\bigoplus\limits_{m_i} \Lambda _{L_i;m_i}
\qquad \text{(resp.} \quad 
 \Lambda _R =\txt\bigoplus\limits^q_{i=1}\bigoplus\limits_{m_i}\Lambda _{R_i;m_i}\ )\]
where $\Lambda _{L_i;m_i}$ (resp. $\Lambda _{R_i;m_i}$~) is the $i$-th \lr Hecke sublattice having multiplicity $\sup(m_i)$~.\\
Indeed, on each \lr weight $\lambda _{L_i}$ (resp. $\lambda _{R_i}$~), which is a character of $\Rep(T_{2j}(L_{v_i}))$ (resp. $\Rep(T^t_{2j}(L_{\o v_i}))$~), there is the action of the Weyl semigroup $W_L$ (resp. $W_R$~) given by:
\[ \phi (s_{i_L}) =w_{i_L} \lambda _{L_i}
\qquad \text{(resp.} \quad 
\phi (s_{i_R})=w_{i_R} \lambda _{R_i}\ )\]
where
\Bi
\item $\phi (s_{i_L})$ (resp. $\phi (s_{i_R})$~) is a \lr Hecke character;
\item $w_{i_L}\in W_L$~, $w_{i_R}\in W_R$~.
\Ei
The \lr action of the Weyl group $W_L$ (resp. $W_R$~) consists in generating the multiplicities 
of the Hecke sublattices $\Lambda _{L_i}$ (resp. $\Lambda _{R_i}$~) to which correspond the subsemimodules $M_{v_i,m_i}$ (resp. $M_{\o v_i,m_i}$~).\epr\Ee
 \vskip 11pt

\begin{coro} \quad The general bilinear Eisenstein cohomology is characterized by the\linebreak K\"unneth isomorphism:
\[ H^E_{R\times_{(D)}L} : \quad
H^{2j}(\partial \o S_{G_R},\widetilde M^{2j}_R)\times_{(D)}H^{2j}(\partial \o S_{G_L},\widetilde M^{2j}_L)
\overset{\sim}{\To} 
 H^{2j}(\partial \o S_{G\RL},\widetilde M^{2j}_R\otimes_{(D)}\widetilde M^{2j}_L) 
\;.\]
\end{coro} \vskip 11pt

\paragraph{Sketch of proof:} \quad this is equivalent to defining the diagonal or complete product between a right and a left linear Eisenstein cohomology semigroup.\epr
 \vskip 11pt

\begin{defi} \quad {\bbf Complete reducibility of $GL_{2n}(L_{\o v}\times L_v)$ \cite{29}:\/} \quad Let 
\begin{align*}
n_L &= 1_{1_L}+\cdots 1_{k_L} + \cdots + 1_{\ell_L}+\cdots+ 1_{n_L}\\
\text{(resp.} \quad
n_R &= 1_{1_R}+\cdots 1_{k_R} + \cdots + 1_{\ell_R}+\cdots+ 1_{n_R}\ )\end{align*}
be a \lr partition of $n_L$ (resp. $n_R$~) labeling the irreducible representations of $T_{2n_L}(L_v)$ (resp. $T_{2n_R}(L_{\o v})$~).\\
Then,
\Be
\item $\Rep(GL_{2n=2_1+\cdots+2_\ell+\cdots+2_n}(L_{\o v}\times L_v)) = \mathop{\boxplus}\limits^{2n}_{2_\ell=2} \Rep(GL_{2_\ell}(L_{\o v}\times L_v))$ 

constitutes a completely reducible {\bf orthogonal bilinear representation\/} of\linebreak $GL_{2n}(L_{\o v}\times L_v)$~;
\item $\Rep(GL_{2n\RL}(L_{\o v}\times L_v))$

$ \qquad  = \mathop{\boxplus}\limits^{2n}_{2_{\ell_R}=2_{\ell_L}=2} \Rep(GL_{2_{\ell\RL}}(L_{\o v}\times L_v))%
\mathop{\boxplus}\limits_{2_{k_R}\neq 2_{\ell_L}} 
\Rep(T^t_{2_{k_{R}}}(L_{\o v})\times T_{2_{\ell_L}}(L_v))$~,

where $GL_{2_{\ell\RL}}$ is another notation for $G_{2_\ell}$~, constitutes a completely 
reducible {\bf nonorthogonal bilinear representation\/} of $GL_{2n}(L_{\o v}\times L_v)$~.
\Ee
\end{defi}
 \vskip 11pt

\begin{propo} \quad Let $\widetilde M^{2n}_L$ (resp. $\widetilde M^{2n}_R$~) be a $2n$-dimensional semisheaf on the $T_{2n}(L_v)$-semimodule (resp. $T^t_{2n}(L_{\o v})$-semimodule).  \\
Let $\partial\o S^{P_{2n=2_1+\cdots+2_n}}_{G_{2n=2_1+\cdots+2_n}}$ and
$\partial\o S^{P_{2n_R\times 2n_L}}_{G_{2n_R\times 2n_L}}$ denote respectively a completely 
reducible orthogonal and nonorthogonal Shimura bisemivariety instead of $\partial \o S_{G\RL}$~.\\
Then, the $2n$-th bilinear Eisenstein cohomologies decompose into direct sums of completely irreducible orthogonal and nonorthogonal bilinear Eisenstein cohomologies according to:
\Bi
\item $H^{2n}(\partial \o S^{P_{2n=2_1+\cdots+2_n}}_{G_{2n=2_1+\cdots+2_n}},\widetilde M^{2n}_R\otimes_D \widetilde M^{2n}_L)$

$\qquad = \bigoplus\limits_{\ell_R=\ell_L} H^{2_\ell}
(\partial \o S ^{P_{2_{\ell_R,\ell_L}}}_{G_{2_{\ell_R,\ell_L}}}  ,\widetilde M^{2_{\ell_R}}_R\otimes \widetilde M^{2_{\ell_L}}_L)$

$\qquad \simeq \Repsp(GL_{2n=2_1+\cdots+2_\ell+\cdots+2_n}(L_{\o v}\times L_v)$~. \vskip 11pt

\item $H^{2n}(\partial \o S^{P_{2n_R\times 2n_L}}_{G_{2n_R\times 2n_L}},\widetilde M^{2n}_R\otimes \widetilde M^{2n}_L)$

$\qquad =\bigoplus\limits_{\ell_R=\ell_L}H^{2_\ell} 
(\partial \o S ^{P_{2_{\ell_R},2_{\ell_L}}}_{G_{2_{n_R},2_{n_L}}}  ,\widetilde M^{2_{\ell_R}}_R\otimes \widetilde M^{2_{\ell_L}}_L)$

$\qquad \quad \bigoplus\limits_{k_R\neq \ell_L} 
H^{2_{k_R},2_{\ell_L}}(\partial \o S^{P_{2_{k_R}\times 2_{\ell_L}}}_{G_{2_{k_R}\times 2_{\ell_L}}},\widetilde M^{2_{k_R}}_R\otimes \widetilde M^{2_{\ell_L}}_L)$

$\qquad \simeq \Repsp (GL_{2n\RL}(L_{\o v}\times L_v))$
\Ei
where $\widetilde M_L^{2_{\ell_L}}$ is a  semisheaf on the $T_{2_{\ell_L}}(L_v)$-semimodule.
\end{propo} \vskip 11pt

\bpr \Be \item The completely reducible orthogonal and nonorthogonal Shimura bisemivarieties are given respectively by:
\begin{align*}
&\partial \o S^{P_{2n=2_1+\cdots+2_n}}_{G_{2n=2_1+\cdots+2_n}}\\
& \quad =\txt\bigoplus\limits^n_{\ell=1} \partial \o S^{P_{2_{\ell_R,\ell_L}}}
_{G_{2_{\ell_R,\ell_L}}}\\
& \quad =\txt\bigoplus\limits^n_{\ell=1} P_{2_{\ell_R,\ell_L}}
({L_{\o v^1}}\times {L_{v^1}})\setminus GL_{2_{\ell_R,\ell_L}}(L^+_{R_d}\times L^+_{L_d})\big/ 
GL_{2_{\ell_R,\ell_L}}((\ZZ/N\ \ZZ)^2)\end{align*}
and by:
\[\partial \o S^{P_{2n_R\times 2n_L}}_{G_{2n_R\times 2n_L}}
=\txt\bigoplus\limits^n_{\ell_R=\ell_L=1} 
\partial \o S^{P_{2_{\ell_R},2_{\ell_L}}}_{G_{2_{\ell_R},2_{\ell_L}}}
\txt\bigoplus\limits^n_{k_R\neq \ell_L=1} \partial \o S^{P_{2_{k_R}\times 2_{\ell_L}}}_{G_{2_{k_R}\times 2_{\ell_L}}}\;.\]

\item The decomposition of the $2n$-th bilinear Eisenstein cohomology into completely irreducible two-dimensional bilinear Eisenstein cohomologies results from its bijection with $\Rep(GL_{2n=2_1+\cdots+2_n}(L_{\o v}\times L_v))$ or with 
$\Rep(GL_{2n\RL}(L_{\o v}\times L_v))$ according to definition 2.5.

\item Every two-dimensional Eisenstein bilinear cohomology decomposes with respect to the places in $(L_{\o v}\times L_v)$ according to one-dimensional components:
\begin{equation} 
H^{2_\ell}
(\partial \o S^{P_{2_{\ell_R},2_{\ell_L}}}_{G_{2_{\ell_R},2_{\ell_L}}},
\widetilde M^{2_{\ell_R}}_R\otimes \widetilde M^{2_{\ell_L}}_L)
\simeq \txt\bigoplus\limits^q_{i=1}
\bigoplus\limits_{m_i}(\widetilde M^{2_{\ell_R}}_{\o v_i;m_i}\otimes \widetilde M^{2_{\ell_L}}_{v_i;m_i})\;.\tag*{\eop}
\end{equation}
\Ee \vskip 11pt

\begin{defi}\quad {\bf Cuspidal representation in terms of global elliptic semimodules:\/}\quad {The decomposition of the Eisenstein bilinear cohomology into one-dimensional irreducible  components needs a cuspidal automorphic representation in terms of global elliptic bisemimodules.

Assume that $f_L$ is a normalized eigenform (of a Hecke operator), holomorphic in the Poincare upper half plane $H$ in $\cit$~, and defined in $\{\Im(z_L)>0\}$~.  $f_L$~, expanded in formal power series $f_L=\sum\limits^q_{i=1}a_{i_L}q_L^i$~, where $q_L=e^{2\pi iz_L}$~, $z_L\in\cit$~, is  a cusp  form of the space $S_L(N)$ and is an eigenvector of the Hecke operators $T_{q_L}$ for $q \nmid N$ and $U_{q_L}$ for $q \mid N$ where $N$ is a positive integer.  The Fourier coefficients $a_{i_L}$ are eigenvalues $c(i,f_L)$ of Hecke operators such that $c(i,f_L)$ generate the ring of integers $\theta _L$ which leads to consider $S_L(N)$ as a $\theta _L$-algebra.

The coalgebra $S_R(N)$ of $S_L(N)$~, defined in the Poincare lower half plane $H^*$~, is then composed of dual cusp forms $f_R=\sum\limits_{i=1}^q a_{i_R}q^i_R$ with $q_R=e^{-2\pi iz_R}$~, $z_r\in\cit$~,  $z_R =z^*_L$~, which are also eigenfunctions of Hecke operators $T_{q_R}$ for $q\nmid N$ and $U_{q_R}$ for $q \mid N$~.

On the other hand, assume that the semisheaf $ \widetilde M^{1_\ell }_{R,L}$ on the semimodule $M^{1_\ell}_{R,L}$ decomposes according to the 
conjugacy classes  of $M^{1_\ell }_{R,L}$  into a set $s_{R,L}=\Gamma (\widetilde M^{1_\ell}_{R,L})$ of   one-dimensional sections 
$s_{i_{R,L}}$~, $1\le i\le q\le \infty$~.  For each section $s_{i_{R,L}}\in s_{R,L}$~, let $\End(G_{s_{R,L}})$ be the Frobenius endomorphism of the group $G_{s_{R,L}}$ of the $s_{i_{R,L}}$ and let $q_{R,L}\to q^i_{R,L}$ be the corresponding Frobenius substitution with $q^i_{R,L}=e^{\pm 2\pi i(i)x}$~, $x\in\rit$~.

A global  elliptic right (resp. left) $G_{s_{R,L}}$-semimodule $\phi_{R,L}(s_{R,L})$ is a ring homomorphism:\linebreak $\phi _{R,L}:s_{R,L}\to \End(G_{s_{R,L}})$ defined by
\[ \phi _{R,L}(s_{R,L})=\txt\sum\limits_i\sum\limits_{m_i} \phi (s_{i_{R,L}})_{i,m_i} q^i_{R,L}\]
where $\sum\limits_i$ runs over the one-dimensional sections of $\widetilde M^{1_\ell}_{R,L}$ and where $\sum\limits_{m_i}$ runs over the ideals of the decomposition group $D_{i^2}$ of the biplace $\o v_i\times v_i$~.

Then, the space $S_{R,L}(\phi _{R_L})$ of global elliptic $G_{s_{R,L}}$-semimodules $\phi _{R,L}(s_{R,L})$ is included into
 the space $S_{R,L}(N)$ of cusp forms $f_{R,L}:S_{R,L}(\phi _{R,L})\hookrightarrow S_{R,L}(N)$ implying that $f_{R,L}\simeq \phi _{R,L}(s_{R,L})$~.}
\end{defi}
\vskip 11pt

\begin{defi}\quad {\bf The decomposition group:\/} \quad {The ring of endomorphisms acting on the global elliptic $G_{s_{R,L}}$-semimodules included into weight two cusp forms is generated over $\ZZ/N\ \ZZ$ by the Hecke operators $T_{q_{R,L}}$ for $q\nmid N$ and $U_{q_{R,L}}$ for $q\mid N$~.  The coset representatives of $U_{q_L}$ are upper triangular and are given by the integral matrices $\left(\begin{smallmatrix} 1&b_N\\ 0&q_N\end{smallmatrix}\right)$ while the coset representatives of $U_{q_R}$ are lower triangular and are given by the matrices 
$\left(\begin{smallmatrix} 1&0\\ b_N&q_N\end{smallmatrix}\right)$~.  
For a general integer $r=a\cdot d$~, we would have respectively the integral matrices $\left(\begin{smallmatrix} a_N & b_N \\ o & d_N\end{smallmatrix}\right)$ and
$\left(\begin{smallmatrix} a_N & 0\\ b_N & d_N\end{smallmatrix}\right)$ of determinant $r\cdot N^2\simeq a_N\cdot d_N$ (class ``0'' of integers modulo $N$~).  These integral matrices modulo $N$ are defined by considering that $q_N=*\mod N\simeq q\cdot N$ and $b_N=*\mod N$~.
Taking into account that the unipotent radical $u(b)=\left(\begin{smallmatrix} 1&b\\ 0&1\end{smallmatrix}\right)$ and its transposed $u(b)^t= \left(\begin{smallmatrix} 1&0\\ b&1\end{smallmatrix}\right)$ generate $\FF_q$~, the following coset representative
\[ GL_2((\ZZ/N\ \ZZ)^2_{|q^2})=\left[\BM 1&b_N\\ 0&1\EM \BM 1&0\\ b_N&1\EM\right]\BM 1&0\\ 0&q_N^2\EM\]
will be adopted for the tensor product $U_{q_R}\otimes U_{q_L}$ of Hecke operators where $\alpha _{q_N^2}= \left(\begin{smallmatrix} 1&0\\ 0&q_N^2\end{smallmatrix}\right)$ is the split Cartan subgroup and where $D_{q_N^2,b_N}=\left(\begin{smallmatrix} 1&b_N\\ 0&1\end{smallmatrix}\right) \left(\begin{smallmatrix} 1&0\\ b_N&1\end{smallmatrix}\right)$ is the representation of the decomposition group associated with $\alpha _{q_N^2}$~.
Then, $GL_2((\ZZ/N\ \ZZ)^2_{|q^2})$ corresponds to a Gauss decomposition of the class ``~$q_N^2$~''.
}\end{defi}
\vskip 11pt

\begin{propo} \quad The eigenvalues $\lambda _{\pm}(q_N^2,b_N^2)$ of the coset representatives\linebreak $GL_2(\ZZ/N\ \ZZ)^2_{|q^2})$ of $U_{q_R}\times U_{q_L}$ are such that they are coefficients of the global elliptic $G_{s_{R,L}}$-semimodules $\phi _{R,L}(s_{R,L})$~, i.e. $\phi (s_{q_{R,L}})_{q,b}\equiv\lambda _{\pm}(q_N^2,b_N^2)$~.  Then, the one-dimensional components of the global elliptic semimodule $\phi _{R,L}(s_{R,L})$ are one-dimensional semitori localized respectively in the upper and in the lower half space and characterized by radii given by $r(q_N^2,b_N^2)=(\lambda _+(q_N^2,b_N^2)-\lambda _-(q_N^2,b_N^2))/2$~.
\end{propo}
\vskip 11pt

\paragraph{Proof:} The eigenvalues of $GL_2(\ZZ/N\ \ZZ)^2_{|q^2})=\alpha_{q_N^2}\cdot D_{q_N^2,b_N^2}$ are
\[ \lambda _{\pm}(q_N^2,b_N^2)=\F{(1+b_N^2+q_N^2)\pm[(1+b_N^2+q_N^2)^2-4q_N^2]^{\half}}2\]
and verify
\begin{align*}
\Tr (GL_2(\ZZ/N\ \ZZ)^2_{|q^2})) &= 1+b_N^2+q_N^2\;, \\
\mbox{and}\quad \det (GL_2(\ZZ/N\ \ZZ)^2_{|q^2})) &= \lambda _+(q_N^2,b_N^2)\cdot \lambda _-(q_N^2,b_N^2)=q^2_N\;.\tag*{\eop}\end{align*}
\vskip 11pt

\begin{propo} \quad According to the Langlands bilinear global program \cite{29} and proposition~2.6, every two-dimensional Eisenstein bilinear cohomology is in bijection with a global elliptic $G_{s\RL}$-bisemimodule $\phi _R(s_R)\otimes \phi _L(s_L)$~:
\begin{align*}
&H^{2_\ell}
(\partial \o S^{P_{2_{\ell_R,\ell_L}}}
_{G_{2_{\ell_R,\ell_L}}}, \widetilde M^{2_{\ell_R}}_{R_\oplus} \otimes \widetilde M^{2_{\ell_L}}_{L_\oplus})
\simeq \txt\bigoplus\limits_{q=1}^q \bigoplus\limits_{m_i} \L(
\widetilde M_{\o v_{i,m_i}}^{1_{\ell_R}}\otimes
\widetilde M_{v_{i,m_i}}^{1_{\ell_L}}\R)
\approx \phi _R(s_R)\otimes \phi _L(s_L)\\
& \quad = \txt \sum\limits^q_{i=1} \sum\limits_{m_i} \lambda _+(i^2_N,m^2_i)e^{-2\pi i(i)x} \otimes \sum\limits^q_{i=1} \sum\limits_{m_i} \lambda _-(i^2_N,m^2_i)e^{2\pi i(i)x} \end{align*}
in such a way that the $i$-th bisection on the $GL_{2_\ell}(L_{\o v_i}\times L_{v_i})$-subbisemimodule 
$(M^{2{\ell_R}}_{R}\otimes M^{2_{\ell_L}}_{L})$ in $H^{2\ell} 
(\partial \o S^{P_{2_{\ell_R},2_{\ell_L}}} _{G_{2_{\ell_R},2_{\ell_L}}},\widetilde M_R^{2_{\ell_R}}\otimes 
\widetilde M_L^{2_{\ell_L}})$ be in one-to-one correspondence with the set of  $m_i$ biterms $\{\lambda _+(i^2_N,m^2_i)e^{-2\pi i(i)x} \times \lambda _-(i^2_N,m^2_i)e^{2\pi i(i)x}\}_{m_i} $ of the global elliptic bisemimodule $\phi _R(s_R)\otimes \phi _L(s_L)$~. So, the global elliptic bisemimodule constitutes the cuspidal representation of the Eisenstein bilinear  cohomology.
\end{propo} \vskip 11pt

\begin{defi}\quad {\bf The semialgebra of elliptic operators ${\rm Ell}_{R,L}(\widetilde 
M_{R,L})$}   is the semialgebra of linear differential operators $D_{R,L}$ defined on
the space $\Gamma_{R,L}(\widetilde M_{R,L})$ of smooth sections $s _{R,L}$ of
$\widetilde M_{R,L}$ and having their principal symbol $\sigma (D_{R,L})$ inversible
\cite{3}.
\end{defi}
\vskip 11pt

\begin{propo}\quad   The bilinear Hilbert semispace $H^{\mp}_a$ is the natural
representation space for the semialgebra of elliptic operators.
\end{propo}
\vskip 11pt

\paragraph{Proof:}  Taking into account the $B_L\circ p_L$ (resp.
$B_R\circ p_R$~) map as introduced in definitions 1.8, the bisemisheaf $\widetilde M_R\otimes_{(D)}\widetilde M_L$ on the 
$GL_n(L_{\o v}\times L_v)$-bisemimodule $M_R\otimes_{(D)} M_L$ will be
transformed into an extended internal or internal left (resp. right)
Hilbert bisemisheaf according to:
\nobeqn B_L\circ p_L &:\quad & \widetilde M_R\otimes_{(D)}\widetilde M_L\To
\widetilde M_{L_R}\otimes_{(D)}\widetilde M_L\equiv \widetilde \Ms_L\;, \\
\noalign{\vskip 6pt}
B_R\circ p_R &:\quad & \widetilde M_R\otimes_{(D)}\widetilde M_L\To
\widetilde M_{R_L}\otimes_{(D)}\widetilde M_R\equiv \widetilde \Ms_R\;.\noeeqn

Consequently, ${\rm Ell}_R(\widetilde M_R)\otimes_{(D)}{\rm
Ell}_L(\widetilde M_L)$ acting on $\widetilde M_{L_R}\otimes_{(D)}\widetilde M_L$ (resp.
$\widetilde M_{R_L}\otimes_{(D)}\widetilde M_R$~) will be an algebra of bioperators (or a
bisemialgebra of operators) acting on an extended internal or internal left
(resp. right) bilinear Hilbert semispace $H^+_a $ (resp. $H^-_a$~) or
$\Hs^+_a$ (resp. $\Hs^-_a$~) and will be noted: \bt[t]{ll}&$[{\rm Ell}_R(\widetilde M_R)\otimes {\rm Ell}_L(\widetilde M_L)](H^{\pm}_a)$\\ or &
$[{\rm Ell}_R(\widetilde M_R)\otimes_D {\rm Ell}_L(\widetilde M_L)](\Hs^{\pm}_a)$~.\te

On the other hand, a semialgebra of operators ${\rm Ell}_{R,L}(\widetilde M_{R,L})$ acting on $H^{\mp}_a$ or $\Hs^{\mp}_a$ will be given by
${\rm Ell}_{R,L}(\widetilde M_{R,L})(H^{\mp}_a)$ or
${\rm Ell}_{R,L}(\widetilde M_{R,L})(\Hs^{\mp}_a)$ in such a way that
${\rm Ell}_{R}(\widetilde M_{R})$ (resp. ${\rm Ell}_{L}(\widetilde M_{L})$~) be a semialgebra of right (resp. left) elliptic linear operators acting
on the set of sections of the semisheaf $\widetilde M_R$ (resp. $\widetilde M_L$~) over the
 $G_R(L_{\o v})$ (resp. $G_L(L_v)$~)-semimodule $M_R$ (resp. $M_L$~)
of $H^{\mp}_a$ or $\Hs^{\mp}_a$~, where $G_R$ (resp. $G_L$~) is another notation for $T^t_n$ (resp. $T_n$~).

Taking into account the considerations given about the enveloping algebras
in definition 1.1, it then becomes clear that the extended bilinear
Hilbert semispace $H^{\pm}_a$ is the natural representation space for the
bisemialgebra and the semialgebra of elliptic operators.\epr
\vskip 11pt

\begin{defis}\quad {\bf a)\ \ Semialgebra of bounded operators:}\quad   
If $\Ls^B_{R,L}(\widetilde M_{R,L})$ denotes the semialgebra of right (resp. left)
operators on the semisheaf $\widetilde M_{R,L}$ over the $G_{R,L}(L_{\o v,v})$-semimodule $M_{R,L}$~, then the
semialgebra of right (resp. left) self-adjoint bounded operators $T_{R,L}$ on
$H^{\mp}_a$ or $\Hs^{\mp}_a$ will be given by: $\Ls^B_{R,L}(H^{\mp}_a)$ and
$\Ls^B_{R,L}(\Hs^{\mp}_a)$~, while the bisemialgebra of self-adjoint bounded
operators on $H^{\mp}_a$ and on $\Hs^{\mp}_a$ will be: $(\Ls^B_R\otimes
\Ls^B_L)(H^{\mp}_a)$ and $(\Ls^B_R\otimes_D
\Ls^B_L)(\Hs^{\mp}_a)$ such that the right and left self-adjoint bounded
operators $T_{R,L}\in\Ls^B_{R,L}$ act respectively on the right and left
semisheaves of $H^{\mp}_a$ and $\Hs^{\mp}_a$~.
\vskip 11pt

\noindent {\bf b)}\ \ A {\bf weight\/} on a semialgebra $\Ls^B_{R,L}(H^+_a)$ is given by
the positive bilinear form $(T_Rs_{i_{L_R}},s_{i_L})$ or
$(s_{i_{L_R}},T_Ls_{i_L})$ which is a map from $\Ls^B_{R,L}(\widetilde M_{L_R}\times
\widetilde M_L)$ into $\cit$ for every section $s _{i_{L_R}}\in \widetilde M_{L_R}$ and
$s _{i_L}\in \widetilde M_L$~.

Similarly, a weight on a bisemialgebra $(\Ls^B_R\otimes \Ls^B_L)(H^+_a)$ will be
given by the positive bilinear form $(T_Rs_{i_{L_R}},T_Ls_{i_L})$ which is a
map from $(\Ls^B_R(\widetilde M_{L_R})\times \Ls^B_L(\widetilde M_L))$ into $\cit$ for all
$T_{R,L}\in \Ls^B_{R,L}$~.
\end{defis}
\vskip 11pt

\begin{defi}\quad {\bf Complex analytic semivariety:}\quad    Let 
$\o X_{S_R}$
(resp. $\o X_{S_L}$~) denote the  \rl complex semispace compactified from
$X_{S_R}=GL_n(L_R^{(\cit)})/GL_n(\ZZ/N\ \ZZ)$
(resp. $X_{S_L}=GL_n(L_L^{(\cit)})/GL_n(\ZZ/N\ \ZZ)$~) being the \rl complex (lattice) analytic semivariety introduced in section 2.1 and covered by 
$\partial \o X_{S_R}$ (resp. $\partial \o X_{S_L}$~) or by
$\partial \o S_{G_R}$ (resp. $\partial \o S_{G_L}$~).

Let $\widetilde M^s_{R,L}$ be an analytic semisheaf on $\o X_{S_R}$ (resp. $\o X_{S_L}$~).

Then, the analytic de Rham  cohomology
$H^*(\o X_{S_{R,L}},\widetilde M^s_{R,L})$ can be computed through the analytic de
Rham complex taking into account that:
\end{defi}
\vskip 11pt

\begin{lm} There is an isomorphism between the (algebraic) Eisenstein
cohomology $H^*(\partial \o S_{G_{R,L}},\widetilde M_{R,L})$ and the analytic de Rham cohomology $H^*(\partial \o X_{S_{R,L}},\widetilde M^s_{R,L})$ .
\end{lm}
\vskip 11pt

\paragraph{Proof:}  Indeed, the isomorphism between the following two de
Rham cohomologies of $\Omega^*$-smooth differential forms with respect to
$\partial\o S_{G_{R,L}}$ and $\o X_{S_{R,L}}$ \cite{20}, \cite{8}, \cite{12}:
\[ H^*(\Omega^*_{\partial \o S_{G_{R,L}}})\simeq
H^*(\Omega^*_{\o X_{S_{R,L}}})\]
leads naturally to the following isomorphism:
\vskip 6pt
\hfill $H^*(\partial \o S_{G_{R,L}},\widetilde M_{R,L})\simeq
H^*(\o X_{S_{R,L}},\widetilde M^s_{R,L})\;.$\epr
\vskip 11pt

\begin{defi}\quad {\bf Analytic bilinear Hilbert semispaces:}\quad    From the
complete (resp. diagonal) bilinear tensor \cite{30} product between the right and left analytic
semisheaves $\widetilde M^s_R$ and $\widetilde M^s_L$~, we can construct a left (resp. right) analytic
bisemisheaf $ \widetilde \Ms^s_{L(D)}$ (resp. $\widetilde \Ms^s_{R(D)}$~) of a left (resp. right) analytic
bilinear Hilbert semispace $H^+_{h}$ or $\Hs^+_{h}$ (resp.
$H^-_{h}$ or $\Hs^-_{h}$~) in complete analogy with which
was done in definition 1.8.
\end{defi}
\vskip 11pt

\begin{propo}\quad The analytic bilinear Hilbert semispace $H^{\pm}_{h}$ is the
natural representation space for the (bi)semialgebras of elliptic operators: ${\rm
Ell}_{R,L}(\widetilde M^s_{R,L})$ and $({\rm
Ell}_{R}(\widetilde M^s_{R})\otimes_{(D)} {\rm
Ell}_{L}(\widetilde M^s_{L}))$~.
\end{propo}
\vskip 11pt

\paragraph{Proof:} This results from definitions 2.12 and 1.1.\epr
\vskip 11pt

\begin{defis}\quad {\bf a)\ \ Serre-Swan theorem:}\quad    Let $\widetilde M^{\rm
top}_{R,L}=C(X_{R,L})$ be the semialgebra of continuous functions
on a compact (semi)space $X_{R,L}$~.  We shall denote by ${\rm
VEC}(X_{R,L})$ the category of complex vector bundles over $X_{R,L}$ and
$P(\widetilde M^{\rm top}_{R,L})$ the category of finitely generated projective right
(resp. left) semimodules $P^{\rm top}_{R,L}$ over $\Gamma ({\rm
VEC}(X_{R,L}))$~.

Then, the Serre-Swan theorem asserts that the categories ${\rm
VEC}(X_{R,L})$ and\linebreak $P^{\rm top}_{R,L}(\Gamma ({\rm
VEC}(X_{R,L})))$ are equivalent \cite{34}.
\vskip 11pt

\noindent {\bf b)\ \ The bisemialgebra $C(X_R\times_{(D)}X_L)$~:}\quad  Let $x_{R,L}$
be a right (resp. left) point of the right (resp. left) semialgebra $\widetilde M^{\rm top}_{R,L}$~.

The complete (resp. diagonal) tensor product between the right and left
semialgebras $\widetilde M^{\rm top}_R$ and $\widetilde M^{\rm top}_L$ can be defined
by:
\[ T^{\rm top}_X : \quad \begin{array}[t]{ccc}
\{\widetilde M^{\rm top}_R,\widetilde M^{\rm top}_L\} &\To & \widetilde M^{\rm top}_R\otimes_{(D)}\widetilde M^{\rm
top}_L\;, \\
\{x_R,x_L\} &\To & x_R\times_{(D)}x_L\;,\end{array}\]
so that the bipoint $x_R\times_{(D)}x_L$
 be characterized by a complete (resp. diagonal) signature.

$\widetilde M^{\rm top}_R\otimes_D\widetilde M^{\rm top}_L$
 is then a finitely generated bisemialgebra.
\vskip 11pt

\noindent {\bf c)\ \ Topological bilinear Hilbert semispace:} \quad   By application of the
$(B_L\circ p_L)$ (resp. $B_R\circ p_R$~) linear map, the
bisemialgebra $(\widetilde M^{\rm top}_R\otimes_{(D)}
\widetilde M^{\rm top}_L)$ can be transformed into an extended internal or internal
left (resp. right) topological Hilbert bisemisheaf $\widetilde \Ms^{\rm top}_{L(D)}$
(resp. $\widetilde \Ms^{\rm top}_{R(D)}$~) which becomes an extended internal or
internal left (resp. right) topological bilinear Hilbert semispace
$H^{\mp}_{\rm top}$ or $\Hs^{\mp}_{\rm top}$ if it is endowed with a
complete or a diagonal bilinear form with values in $\cit$~.
\end{defis}
\vskip 11pt

\begin{propo}\quad The extended internal and internal left (resp. right)
topological bilinear Hilbert semispaces $H^{\pm}_{\rm top}$ and
$\Hs^{\pm}_{\rm top}$ are $C^*$-(bi)semialgebras.
\end{propo}
\vskip 11pt

\paragraph{Proof:}  By definition $\widetilde M^{\rm top}_{R,L}$ is a right (resp.
left)   semialgebra $C(X_{R,L})$ of
continuous sections $s^{\rm top}_{R,L}(X_{R,L})$ on $X_{R,L}$~.

Now, the bisemialgebra $\widetilde \Ms^{\rm top}_{L(D)}$ or
$\widetilde \Ms^{\rm top}_{R(D)}$ is an involutive  bisemialgebra over
$\cit$ of continuous bifunctions $s^{\rm top}_R(X_R)\otimes_{(D)}s^{\rm
top}_L(X_L)$~.  Indeed, the involution, which must be taken into account, is
a bilinear map transforming $H^+_{\rm top}$ or $\Hs^+_{\rm top}$ (resp.
$H^-_{\rm top}$ or $\Hs^-_{\rm top}$~) into $H^-_{\rm top}$ or $\Hs^-_{\rm top}$
(resp. $H^+_{\rm top}$ or $\Hs^+_{\rm top}$~).

Recall the composition of maps:
\nobeqn B_L\circ p_L &:\quad & \widetilde M^{\rm top}_R\otimes_{(D)}\widetilde M^{\rm top}_L\To
\widetilde M^{\rm top}_{L_R}\otimes_{(D)}\widetilde M^{\rm top}_L\equiv \widetilde \Ms^{\rm top}_{L(D)}\;, \\
\noalign{\vskip 6pt}
B_R\circ p_R &:\quad& \widetilde M^{\rm top}_R\otimes_{(D)} \widetilde M^{\rm top}_L\To
\widetilde M^{\rm top}_{R_L}\otimes_{(D)}\widetilde M^{\rm top}_R\equiv \widetilde\Ms^{\rm top}_{R(D)}\;,\noeeqn
as introduced in definitions 1.8.

So, the bilinear map:
\nobeqn (p^{-1}_L\circ B_L^{-1} )\otimes_{(D)}(B_R\circ p_R  ) &:\quad & \widetilde\Ms^{\rm top}_{L(D)}\To {\widetilde\Ms}^{\rm top}_{R(D)}\;, \\
\noalign{\vskip 6pt}
(B_R\circ p_R )\otimes_{(D)}(p^{-1}_L\circ B_L ^{-1} ) &:& \widetilde\Ms^{\rm
top}_{R(D)}\To  {\widetilde\Ms}^{\rm top}_{L(D)}\;, \noeeqn
transforms the bisemialgebra $\widetilde\Ms^{\rm top}_{L(D)}$ (resp. $\widetilde\Ms^{\rm
top}_{R(D)}$~) into the bisemialgebra $\widetilde\Ms^{\rm top}_{R(D)}$ (resp. 
${\widetilde\Ms}^{\rm top}_{L(D)}$~) which corresponds to an antilinear involution
transforming the left (resp. right) bilinear Hilbert semispace $H^+_{\rm top}$
or $\Hs^+_{\rm top}$ (resp. $H^-_{\rm top}$
or $\Hs^-_{\rm top}$~) into the right (resp. left) involuted bilinear
Hilbert semispace $H^-_{\rm top}$
or $\Hs^-_{\rm top}$ (resp. $H^+_{\rm top}$
or $\Hs^+_{\rm top}$~).\epr
\vskip 11pt

\begin{defis}\quad {\bf a)\ \ $K$-functor of Kasparov \cite{25}:}\quad     We are now
interested in extensions of the bisemialgebra $\widetilde M^{\rm
top}_{R,L}$~.  Let $\Ls^B_{R,L}(\widetilde M^{\rm top}_{R,L})$ denote the semialgebra of
bounded operators on $\widetilde M^{\rm top}_{R,L}$ and let $\Ks_{R,L}$ be the ideal of
compact operators.

The set of extension classes of $\Ks_{R,L}$ by $\Ls^B_{R,L}(\widetilde M^{\rm
top}_{R,L})$~, noted ${\rm Ext}(\Ls^B_{R,L}(\widetilde M^{\rm
top}_{R,L}),\Ks_{R,L})$~, is an abelian semigroup naturally isomorphic to
${\rm Ext}(X_{R,L})$ as developed by Brown, Douglas and Fillmore \cite{10},
\cite{11}.

In connection with the work of Atiyah \cite{4}, \cite{5}, G.G. Kasparov constructed a
general $K$-functor $K_*K^*(\widetilde M^{\rm top}_{R,L},\Ls^B_{R,L})$~, special
cases of which are the ordinary cohomological $K$-functor $K^*(\widetilde M^{\rm
top}_{R,L})$ and the homological $K$-functor $K_*(\Ls^B_{R,L})$~.

Especially interesting is the case where the $C^*$-semialgebras $\widetilde M^{\rm
top}_{R,L}$ and $\Ls^B_{R,L}$ are equipped with the continuous action of a
locally compact semigroup $G^{\ell}_{R,L}$~.  This allows to define an
abelian group $KK^{G^\ell_{R,L}}(\widetilde M^{\rm top}_{R,L},\Ls^B_{R,L})$ \cite{26}.
\vskip 11pt

\noindent {\bf b)\ \ Bisemialgebra of bounded operators:}\quad   Considering the
$C^*$-bisemialgebra $\widetilde\Ms^{\rm top}_{L(D)}$ (resp.
$\widetilde\Ms^{\rm top}_{R(D)}$~), the bisemialgebra of bounded operators on it will be
$(\Ls^B_R\otimes_{(D)}\Ls^B_L)(\widetilde\Ms^{\rm top}_{L(D)})$ (resp.
$(\Ls^B_R\otimes_{(D)}\Ls^B_L)(\widetilde\Ms^{\rm top}_{R(D)})$~) or
$(\Ls^B_R\otimes_{(D)}\Ls^B_L)(H_{\rm top}^{\mp})$
(resp. $(\Ls^B_R\otimes_{(D)}\Ls^B_L)(\Hs_{\rm top}^{\mp})$~) if we
envisage their actions on the extended (resp. diagonal) bilinear Hilbert
semispace $H^{\mp}_{\rm top}$ (resp. $\Hs^{\mp}_{\rm top}$~).
\end{defis}
\vskip 11pt

\section{Von Neumann semialgebras and bisemialgebras}

\begin{defis}\quad {\bf a)\ \ Norm topology of bounded operators:}\quad    Let
$(\Ls^B_R\otimes\Ls^B_L)(H^{+}_{\rm top})$ be the bisemialgebra of bounded
operators acting from the topological extended bilinear Hilbert semispace
$H^+_{\rm top}$ into itself.

Then, the norm topology for an operator $T_R\otimes T_L\in \Ls^B_R\otimes
\Ls^B_L$ will be defined by
\[\|T_R\otimes T_L\|=
\sup \L( {\|T_Rs ^{\rm top}_{L_R}\times T_Ls ^{\rm
top}_L\|\big/ \| s ^{\rm top}_{L_R}\times s ^{\rm
top}_L \|}\R)\;,\]
for every section $s ^{\rm top}_{L_R}\in \widetilde M^{\rm top}_{L_R}$ and 
$s ^{\rm top}_{L}\in \widetilde M^{\rm top}_{L}\subset H^+_{\rm top}$~, since, if $\Ls^B_L(H^+_{\rm top})$ is 
the semialgebra of left bounded operators acting on the semisheaf $\widetilde M^{\rm
top}_L$ of $H^+_{\rm top}$~, the norm topology for a left bounded
operator $T_L$ is given by
\[\|T_L\| =\sup\L(  {\|T_Ls ^{\rm
top}_L\|\big/ \|s ^{\rm
top}_L\|}\R) \;.\]
\vskip 11pt

\noindent {\bf b)}\ \ An {\bf involution\/} on the operator $T_{R,L}$ is defined by
\nobeqn i_L &:\quad & T_R\To T^\dagger_R\equiv T_L\;, \\
\noalign{\vskip 6pt}
 i_R &:\quad & T_L\To T^\dagger_L\equiv T_R\;, \noeeqn
such that $(T^\dagger_Rs ^{\rm top}_{L_R},T^\dagger_Ls ^{\rm top}_L)=
(T_Rs ^{\rm top}_{L_R},T_Ls^{\rm top}_L)$ making $T_R$ and $T_L$
self-adjoint.
\end{defis}
\vskip 11pt

\begin{defis} $\;$ {\bf Bisemialgebras of von Neumann on extended bilinear Hilbert
semispaces:}
\quad {\bf a)\/}\ A right (resp. left) semialgebra of von Neumann
$\MM_{R,L}(H^{\mp}_{\rm top})$ in the  topological
extended bilinear Hilbert semispace $H^{\mp}_{\rm top}$ is an involutive
subalgebra of $\Ls^B_{R,L}(H^{\mp}_{\rm top})$ having a closed norm
topology \cite{19}.

Similarly, a semialgebra of von Neumann $\MM_{R,L}(H^{\mp}_{h})$ in $H^{\mp}_{h}$ is an involutive subsemialgebra of ${\rm
Ell}_{R,L}(\widetilde M^s_{R,L})$ having a closed norm topology.
\vskip 11pt

\noindent {\bf b)\/}\ \ A bisemialgebra of von Neumann 
$\MM_{R\times L}(H^{\mp}_{\left\{{\rm top}\atop{h}\right.})$
 in $H^{\mp}_{\left\{{\rm top}\atop{h}\right.})$ is an involutive subbisemialgebra of 
$(\Ls^B_R\otimes
\Ls^B_L)(H^{\mp}_{\left\{{\rm top}\atop{h}\right.})$ having a closed norm topology \cite{19}.
\vskip 11pt

\noindent {\bf c)\/}\ \ A bisemialgebra of von Neumann $\MM_{R\times L}(H^{\mp}_a)$ in the
algebraic extended bilinear Hilbert semispace $H^{\mp}_a$ is also an
involutive subbisemialgebra of $(\Ls^B_R\otimes \Ls^B_L)(H^{\mp}_a)$ having a
closed norm topology.
\end{defis}
\vskip 11pt

\begin{propo}\quad  Let $\MM_{R,L}(H^{\mp}_a)$ and
$\MM_{R\times L}(H^{\mp}_a)$ be respectively a semialgebra and a bisemialgebra
of von Neumann on the algebraic bilinear extended Hilbert semispace
$H^{\mp}_a $~.

Let $\MM_{R,L}(H^{\mp}_{h})$ and
$\MM_{R\times L}(H^{\mp}_{h})$ be respectively a
semialgebra and a bisemialgebra of von Neumann on the analytic bilinear extended
Hilbert semispace $H^{\mp}_{h}$~.

Then, we have the isomorphisms:
\nobeqn i_{\MM^a_{R,L}-\MM^{h}_{R,L}} &:\quad & \MM_{R,L}(H^{\mp}_a)\To
\MM_{R,L}(H^{\mp}_{h})\;, \\
\noalign{\vskip 6pt}
i_{\MM^a_{R\times L}-\MM^{h}_{R\times L}} &:\quad &
\MM_{R\times L}(H^{\mp}_a)\To \MM_{R\times L}(H^{\mp}_{h})\;.\noeeqn
\end{propo}
\vskip 11pt

\paragraph{Proof:} This results immediately from lemma 2.15.\epr
\vskip 11pt

\begin{propo}\quad There exists an isomorphism
\[ i_{M_{\left\{{\rm top}\atop{h}\right.}-\MM^{\left\{{\rm top}\atop{h}\right.}_{R,L}}:\quad M_{\left\{{\rm top}\atop{h}\right.}\To \MM^{\left\{{\rm top}\atop{h}\right.}_{R,L}(H^{\mp}_{\left\{{\rm top}\atop{h}\right.})\]
between an algebra of von Neumann $M_{\left\{{\rm top}\atop{h}\right.}$ on a linear Hilbert 
space ${\it h}_{\left\{{\rm top}\atop{h}\right.}$ {\rm \cite{26}\/} and a semialgebra of von Neumann $\MM^{\left\{{\rm top}\atop{h}\right.}_{R,L}(H^{\mp}_{\left\{{\rm top}\atop{h}\right.})$ on the extended bilinear Hilbert semispace $H^{\mp}_{\left\{{\rm top}\atop{h}\right.}$~.
\end{propo}
\vskip 11pt

\paragraph{Proof:} Let $V_{\left\{{\rm top}\atop{h}\right.}$ be a compact manifold of class $C^0$ (resp. $C^\infty $~)
associated with $M_R\otimes_{(D)} M_L$ in the sense of \cite{29}
and let $A_{\left\{{\rm top}\atop{h}\right.}$ be the corresponding stellar algebra of $C^0$ (resp.
$C^\infty $~) functions on $V_{\left\{{\rm top}\atop{h}\right.}$ with values in $\cit$~.

Then, a Fredholm module on $A_{\left\{{\rm top}\atop{h}\right.}$ is essentially given by the involutive representation $\Pi_{\left\{{\rm top}\atop{h}\right.}$ of $A_{\left\{{\rm top}\atop{h}\right.}$ in a linear Hilbert space $h_{\left\{{\rm top}\atop{h}\right.}$ and by a self-adjoint operator $F$~.

Furthermore, an algebra of von Neumann $M_{\left\{{\rm top}\atop{h}\right.}$ in a linear Hilbert space $h_{\left\{{\rm top}\atop{h}\right.}$ 
is an involutive subalgebra of bounded operators $\Ls(h_{\left\{{\rm top}\atop{h}\right.})$
from $h_{\left\{{\rm top}\atop{h}\right.}$ to $h_{\left\{{\rm top}\atop{h}\right.}$ such that $M_{\left\{{\rm top}\atop{h}\right.}$ be $\sigma (\Ls(h_{\left\{{\rm top}\atop{h}\right.}),\Ls(h_{\left\{{\rm top}\atop{h}\right.})_*)$ closed.

Now, it is clear that there is a one-to-one correspondence between:
\Bean\item a Fredholm module on $A_{\left\{{\rm top}\atop{h}\right.}$ and a subsemialgebra of $\Ls^B_{R,L}(H^{\mp}_{\left\{{\rm top}\atop{h}\right.})$ since the extended bilinear Hilbert semispace $H^{\mp}_{\left\{{\rm top}\atop{h}\right.}$ can be considered as a representation space of the linear
Hilbert space $h_{\left\{{\rm top}\atop{h}\right.}$ covered by $H^{\mp}_{\L\{{{\rm top}\atop{h}}\R.}$ \cite{29}.
\item the weak topological condition of closeness of $\sigma (\Ls(h_{\left\{{\rm top}\atop{h}\right.}),\Ls(h_{\left\{{\rm top}\atop{h}\right.})_*)$ and the condition of closed norm topology of
$\Ls^B_{R,L}(H^{\mp}_{\left\{{\rm top}\atop{h}\right.})$ since $\Ls(h_{\left\{{\rm top}\atop{h}\right.})_*$ is the dual of $\Ls(h_{\left\{{\rm top}\atop{h}\right.})$~.\Ee

As there is an isomorphism between a Fredholm module on $A_{\left\{{\rm top}\atop{h}\right.}$ 
and a subspace of $\Ls(h_{\left\{{\rm top}\atop{h}\right.})$~, we have the announced isomorphism 
$i_{M_{\left\{{\rm top}\atop{h}\right.}-\MM^{\left\{{\rm top}\atop{h}\right.}_{R,L}}\!:\!M_{\left\{{\rm top}\atop{h}\right.}\!\to\! \MM^{\left\{{\rm top}\atop{h}\right.}_{R,L}(H^{\mp}_{\left\{{\rm top}\atop{h}\right.})$~.\epr
\vskip 11pt

\sub{Shifted actions of differential bioperators on the representation spaces of bilinear semigroups}  \Be\item Let $T^{D_m}_{R,L}\in
\MM_{R,L}(H^{\mp}_a)$ be a \rl differential linear operator of rank $m$ (i.e. operating on 
$m$ variables) of the semialgebra of von Neumann $\MM_{R,L}(H^{\pm}_a)$~. This operator 
$T^{D_m}_R$ (resp. $T^{D_m}_L$~), noted in condensed form $T^D_R$ (resp $T^D_L$~), is assumed to 
be associated with the action of a $T^t_m(\rit)$-semigroup structure (resp. a $T_m(\rit)$-semigroup structure) on the \rl $n$-dimensional semisheaf $\widetilde M_R$ (resp. $\widetilde M_L$~) of the extended bilinear Hilbert semispace $H^+_a$~, with $m\le n$~.  Recall \cite{3} that a $T^t_m(\rit)$-semigroup structure (resp. a $T_m(\rit)$-semigroup structure) on $\widetilde M_R$ (resp. $\widetilde M_L$~) means a principal $T^t_m(\rit)$-bundle (resp. a $T_m(\rit)$-bundle) on $\widetilde M_R$ (resp. $\widetilde M_L$~).

\item Similarly, $(T^D_R\otimes T^D_L)$ will denote the tensor product of the right and left 
differential operators (~$T^D_R$ and $T^D_L$~) acting on the  bisemisheaf $(\widetilde M_R\otimes \widetilde M_L)$
 such that $(T^D_R\otimes T^D_L)\in \MM\RL(H^{\pm}_a)$ be associated with a principal $GL_m(\rit\times \rit)=T_m(\rit)\times T^t_m(\rit)$-bundle on $(\widetilde M_R\otimes \widetilde M_L)$~.

%%%%%%%%%%%%%%§ 3.7 avancé ici
%\sub{Shifted actions of bioperators on the representations of bilinear semigroups}
%\Bena

\item Let $(T^D_R\times T^D_L)$ be the tensor product of a right and a left linear differential operator of rank $m$ such that the action of $(T^D_R\otimes T^D_L)$ be associated with a $GL_m(\rit\times\rit)$-principal bundle on the bisemisheaf
$ (\widetilde M_R\otimes  \widetilde M_L)$ over the 
$GL_n(L_{\o v}\times L_v)$-bisemimodule $(M_R\otimes M_L)$~.  Then, the action of $(T^D_R\otimes T^D_L)$ on $( \widetilde M_R\otimes \widetilde M_L)$ is equivalent to:
\Be
\item consider the mapping
\[ T^D_R\otimes T^D_L : \quad \widetilde M_R\otimes \widetilde M_L \To \widetilde M_{R_{n[m]}}\otimes \widetilde M_{L_{n[m]}}\]
from the 
 bisemisheaf
$ (\widetilde M_R\otimes  \widetilde M_L)$ over the
$GL_n(L_{\o v}\times L_v)$-bisemimodule $(M_R\otimes M_L)$ to the 
 bisemisheaf
$ (\widetilde M_{R_{n[m]}}\otimes  \widetilde M_{L_{n[m]}})$ over the 
$GL_{n[m]}((L_{\o v}\otimes \rit)\times (L_{v}\otimes \rit))$-bisemimodule $ (M_{R_{n[m]}}\otimes M_{L_{n[m]}})$ such that 
$(\widetilde M_{R_{n[m]}}\otimes \widetilde M_{L_{n[m]}})$ be a  bisemisheaf shifted into $(m\times m)$ dimensions.

\item consider a shift into $(m\times m)$ dimensions of the functional representation space
$\FRepsp(GL_n(L_{\o v}\times L_v))$ of the general bilinear semigroup $GL_n(L_{\o v}\times L_v)$ 
leading to the homomorphism:
\[T^D_R\otimes T^D_L : \quad \FRepsp(GL_n(L_{\o v}\times L_v)) \To \FRepsp(GL_{n[m]}((L_{\o v}\otimes \rit)\times (L_{v}\otimes \rit)))\]
where $\FRepsp(GL_{n[m]}((L_{\o v}\otimes \rit)\times (L_{v}\otimes \rit)))$~, denoting the functional representation space of $GL_n(L_{\o v}\times L_v)$ shifted into $(m\times m)$ dimensions, is the shifted bisemisheaf  $(\widetilde M_{R_{n[m]}}\otimes  \widetilde M_{L_{n[m]}})$ on the bisemigroupoid $GL_{n[m]} ((L_{\o v}\otimes \rit)\times (L_v\otimes\rit))$ and is equal to:
\begin{multline*}
 \FRepsp(GL_{n[m]}((L_{\o v}\otimes \rit)\times (L_{v}\otimes \rit)))\\
= \AdFRepsp(GL_{m}(\rit\times\rit))\times 
\FRepsp(GL_{n}(L_{\o v}\otimes L_v)\end{multline*}
in such a way that \cite{xxx}
\Bi
\item $\AdFRepsp(GL_{m}(\rit\times\rit))$~, being the adjoint functional representation space of $GL_m(\rit\times\rit)$~, corresponds to the action of $(T^D_R\times T^D_L)$~;
\item $\FRepsp(GL_{n}(L_{\o v}\otimes L_v))$~, being the functional representation space of\linebreak $GL_n(L_{\o v}\times L_v)$~, correspond to the bisemisheaf $(\widetilde M_R\otimes \widetilde M_L)$~.
\Ei
\Ee

\item Similarly, the shifting ``action'' of $(T^D_R\otimes T^D_L)$ on 
functional representation space of
the bilinear subsemigroup $GL_n((\ZZ/N\ \ZZ)^2)$ would be:
\begin{align*}
&T^D_R\otimes T^D_L : \\ 
&\quad \FRepsp ( GL_n(( \ZZ/N\ \ZZ )^2 ))  \\
& \qquad
=\FRepsp ( D_n (( \ZZ/N\ \ZZ )^2 )\times [UT_n^t ( \ZZ/N\ \ZZ )\times UT_n ( \ZZ/N\ \ZZ)])\\
& \qquad \quad\To
\FRepsp(GL_{n[m]} ((\ZZ/N\ \ZZ)^2\otimes\rit^2))\\
&\hspace{2.5cm} = \FRepsp(D_{n[m]}((\ZZ/N\ \ZZ)^2\otimes\rit^2)\\
&\hspace{4cm}\times [UT_{n[m]}^t((\ZZ/N\ \ZZ)\otimes\rit^2)\times UT_{n[m]}((\ZZ/N\ \ZZ)\otimes\rit^2))]\end{align*}
where:
\Bi
\item $\FRepsp(D_{n[m]}((\ZZ/N\ \ZZ)^2\otimes\rit^2))$ is the functional representation space of the subgroup of integer diagonal matrices of 
order $n$ shifted into $m$ dimensions
.
\item $\FRepsp(UT_{n[m]}((\ZZ/N\ \ZZ)\otimes\rit))$ is the functional representation space of the subgroup of integer unitriangular matrices shifted in $m$ dimensions.
\Ei

\item And, the functional representation space of bilinear parabolic subsemigroup $P_n({L_{\o v^1}}\times {L_{v^1}})$ would also be shifted into $(m\times m)$ dimensions under the action of $(T^D_R\otimes T^D_L)$ according to:
\begin{align*}
&T^D_R\otimes T^D_L : \\ 
&\quad  \FRepsp (P_n(({L _{\o v^1}} \times {L _{v^1}})) 
= \FRepsp(D_n( {L _{\o v^1}} \times {L_{v^1}})  
\times [UT_n^t( {L _{\o v^1}} )\times UT_n( {L _{v^1}}  )])\\
& \quad \To
\FRepsp( P_{n[m]}(({L _{\o v^1}}\otimes \rit) \times ({L _{v^1}}\otimes \rit))) \\
&\hspace{25mm} = \FRepsp(D_{n[m]}( ({L _{\o v^1}} \otimes \rit)\times ({L_{v^1}}\otimes \rit)) \\
&\hspace{5cm} 
\times [UT_{n[m]}^t( {L _{\o v^1}} \otimes \rit)\times UT_{n[m]}( {L _{v^1}} \otimes \rit) )])\;.
\end{align*}
%\Ee \vskip 11pt

%%%jusqu'ici

\item On the other hand, referring to section 1.5, $GL_n(L^{nr}_{\o v}\times L^{nr}_v)$ has $GL(M^{nr}_R\otimes M^{nr}_L)\equiv\Gamma _R\times  \Gamma _L$ for bilinear 
(semi)\-group of  automorphisms and has for pseudounramified conjugacy classes the biclasses 
$\gamma {(i)}_R\times  \gamma {(i)}_L$~, $1\le i\le q$~, if the set of fixed bielements  is the smallest normal bilinear subsemigroup $P_n(L^{nr}_{\o v^1}\times L^{nr}_{v^1})$  of $M^{nr}_R\otimes M^{nr}_L$~.  This smallest normal bilinear subsemigroup of $GL_n(L^{nr}_{\o v}\times L^{nr}_v)$ is the $n$-dimensional equivalent of the product, right by left, of the global inertia subgroups ${\rm I}_{L_{\o v_i}}\times {\rm I}_{L_{v_i}}$ of degree $N^2=1$ as introduced in definition 1.4 \cite{27}.

%%%%%%%%petite rajoute
In this context, the action of $(T^D_R\otimes T^D_L)$ on $(\widetilde M^{nr}_R\otimes\widetilde M^{nr}_L)$~, associated with the principal 
$GL_m(\rit\times\rit)$-bundle on $(\widetilde M^{nr}_R\otimes\widetilde M^{nr}_L)$ with group
$GL_m(\rit\times\rit)$~, leads to envisage that the bilinear semigroupoid $GL_{n[m]}((L^{nr}_{\o v}\otimes \rit)\times (L^{nr}_v\otimes\rit))$~, shifting in $(m\times m)$ dimensions, has $GL(\widetilde M^{nr}_{R_{n[m]}}\otimes \widetilde M^{nr}_{L_{n[m]}})
\equiv \Gamma ^{[m]}_R\times \Gamma ^{[m]}_L$ for bilinear semigroupoid of shifted  automorphisms and has for pseudounramified conjugacy classes the biclasses $(\gamma ^{[m]}(i)_R\times \gamma ^{[m]}(i)_L)$ shifted in $(m\times m)$ dimensions, if the set of shifted fixed bielements corresponds to the smallest normal bilinear subsemigroupoid $P_{n[m]}((L^{nr}_{\o v^1}\otimes\rit)\times (L^{nr}_{v^1}\otimes\rit))$~, i.e. the bilinear pseudounramified parabolic subsemigroupoid.

The shifted pseudounramified conjugacy biclasses $(\gamma ^{[m]}(i)_R\times \gamma ^{[m]}(i)_L)$ are in one-to-one correspondence with their unshifted equivalents 
$(\gamma (i)_R\times \gamma (i)_L)$ because the bilinear subsemigroup 
$(\Gamma ^{[m]}_R\times \Gamma ^{[m]}_L)$ of  automorphisms shifting in $(m\times m)$ real dimensions results from the principal $GL_m(\rit\times\rit)$-bundle on $(\widetilde M^{nr}_R\times \widetilde M^{nr}_L)$ and corresponds to the $(m\times n)$-dimensional representation of the product, right by left, of the differential Galois semigroups of the algebraic extensions $L_R^{nr,+}$ and $L_L^{nr,+}$~.

%%%%%%%%

\item $GL_n(L_{\o v}\times L_v)$ has for bilinear subsemigroup of  automorphisms $\Pra \Gamma _R\times \Pra \Gamma _L$ and has for pseudoramified conjugacy classes the biclasses $g{(i)}_R\times g{(i)}_L$ if the set of fixed bielements is of dimension $N>1$ with respect to the basis of $M_R\otimes M_L$~.  These fixed bielements of $g{(i)}_R\times g{(i)}_L$ correspond to the product, right by left, of completions of degrees equal to $N>1$~.

%%%rajoute :

Similarly, $GL_{n[m]}((L_{\o v}\otimes\rit)\times (L_v\otimes \rit))$ has for bilinear subsemigroup of shifted  automorphisms $(\Pra \Gamma ^{[m]}_R\times \Pra \Gamma _L^{[m]})$ and has for shifted pseudoramified conjugacy classes the biclasses $(g^{[m]}{(i)}_R\times g^{[m]}{(i)}_L)$~.  

As $(\Pra \Gamma ^{[m]}_R\times \Pra \Gamma _L^{[m]})$ is the bilinear subsemigroupoid  $GL(\widetilde M_{R_{n[m]}}\otimes\widetilde M_{L_{n[m]}})$ of automorphisms shifting in $(m\times m)$ real dimensions with respect to the biaction of
$(T^D_R\otimes T^D_L)$ on $(\widetilde M_R\otimes \widetilde M_L)$~, associated with the $GL_m(\rit\times\rit)$-principal bundle introduced in 3.), it is clear that the shifted pseudoramified conjugacy biclasses $(g^{[m]}(i)_R\times g^{[m]}(i)_L)$ are in one-to-one correspondence with the unshifted pseudoramified conjugacy biclasses $g (i)_R\times g(i)_L$~.
%%%%%%%

\Ee \vskip 11pt

\begin{propo}  \quad The action of the differential bioperator $(T^D_R\otimes T^D_L)$ of rank $(m\times m)$~, associated with a principal $GL_m(\rit\times \rit)$-bundle on the $(n\times n)$-dimensional pseudo(un)\-ramified bisemisheaf $(\widetilde M^{(nr)}_{R}\otimes \widetilde M^{(nr)}_{L})$~, consists in mapping $(\widetilde M^{(nr)}_{R_\oplus}\otimes \widetilde M^{(nr)}_{L_\oplus})$ into  $(\widetilde M^{(nr)}_{R_{n[m]_\oplus}}\otimes \widetilde M^{(nr)}_{L_{n[m]_\oplus}})$ shifted into $(m\times m)$ dimensions:
\[ T^D_R\otimes T^D_L: \quad \widetilde M^{(nr)}_{R_\oplus}\otimes \widetilde M^{(nr)}_{L_\oplus} \To (\widetilde M^{(nr)}_{R_{n[m]_\oplus}}\otimes \widetilde M^{(nr)}_{L_{n[m]_\oplus}})\]
such that:
\Bean
\item $\widetilde M_{R_{n[m]_\oplus}}\otimes \widetilde M_{L_{n[m]_\oplus}}$ decomposes into shifted pseudoramified 
subbisemisheaves according to the shifted pseudoramified conjugacy biclasses $g^{[m]}_R(i)\times g^{[m]}{(i)}_L$ of the bisemigroupoid
$GL_{n[m]}((L_{\o v}\otimes\rit)\times (L_v\otimes\rit))$ and with respect to the shifted  automorphisms $\Pra \Gamma ^{[m]}_R\times \Pra \Gamma ^{[m]}_L$ of $GL_{n[m]}((L_{\o v}\otimes\rit)\times (L_v\otimes\rit))$ as follows:
\[\widetilde M_{R_{n[m]_\oplus}}\otimes \widetilde M_{L_{n[m]_\oplus}}=
\txt\bigoplus\limits^q_{i=1}
\txt\bigoplus\limits_{m_i}
(\widetilde M_{R_{n[m]}}(i)\otimes \widetilde M_{L_{n[m]}}(i))\]
where the integer $q$ is related to the dimension $(q\cdot N)^n$ of the algebraic basis of $\widetilde M_{R_{n[m]}}(q)$ and $ \widetilde M _{L_{n[m]}}(q)$~, i.e. to the number of Galois  automorphisms.

\item $\widetilde M^{nr}_{R_{n[m]_\oplus}}\otimes \widetilde M^{nr}_{L_{n[m]_\oplus}}$  decomposes  into shifted pseudounramified  subbisemisheaves according to the shifted pseudounramified conjugacy biclasses $\gamma_R ^{[m]}{(i)}\times \gamma_L ^{[m]}{(i)}$ of $GL_{n[m]}((L^{nr}_{\o v}\otimes\rit)\times (L^{nr}_v\otimes\rit))$ as follows:
\[\widetilde M^{nr}_{R_{n[m]_\oplus}}\otimes \widetilde M^{nr}_{L_{n[m]_\oplus}}=\txt\bigoplus\limits^q_{i=1}
\bigoplus\limits_{m_i} ( \widetilde M^{nr}_{R_{n[m]}}(i)\otimes  \widetilde M^{nr}_{L_{n[m]}}(i))\]
where the integer $q$~, i.e. the global class residue degree $f_{v_q}=q$ (see definition 1.3), refers to the algebraic dimension $q^n$ of 
$\widetilde M^{nr}_{R_{n[m]}}(q)$ and of
$\widetilde M^{nr}_{L_{n[m]}}(q)$~.

%\item $t\ge q$~; $t\equiv q$ if $N=1$~.
\Ee
\end{propo} \vskip 11pt

\bpr \Be\item The shifted  bisemisheaf $ (\widetilde M^{(nr)}_{R_{n[m]}}\otimes  \widetilde M^{(nr)}_{L_{n[m]}})$ 
is a biobject of the derived category $D( \widetilde M_R\otimes  \widetilde M_L,\rit\otimes\rit)$~.

\item The algebraic dimension $(q\cdot N)^n$ of $M_{R_{n[m]}}(q)$ and of $ M_{L_{n[m]}}(q)$ corresponds to the number of Galois  automorphisms while the algebraic dimension $(q\cdot N)^m$ corresponds to the number of shifted automorphisms.

\item The pseudounramified algebraic dimension  $q^n$ is   such  that $q$ corresponds to the number of archimedean places of the semifields $L^+_L$ and $L^+_R$~.\epr 
\Ee\vskip 11pt

%%%%%%%%%%%%
\begin{defi} \quad {\bf Pseudoramified  and pseudounramified algebraic dimensions:\/} \quad Until now, two kinds of algebraic dimensions have emerged:
\Bean
\item the ``pseudoramified'' algebraic dimensions $i^n\ N^n$~, referring to the Galois extension degrees being multiples of $N>1$.

The shifted pseudoramified algebraic dimensions  $i^m\cdot N^m$ referring to the dimensions of the $m$-dimensional representations of the differential Galois subgroups;
\item the pseudounramified algebraic dimension $i^n$ referring to the $n$-th powers of the global residue degree $i$~.

The pseudounramified algebraic dimension $i^m$ referring to the dimensions of the $m$-dimensional representations of the corresponding differential Galois subgroups.
\Ee
Consider for example the left $T_n(L_{v_i})$-subsemimodule $M_{v_i}\subset M_L$ (see section 
1.5) having a rank $n_i=i^n\cdot N^n=f^n_{v_i}\cdot N^n$~.  Then, the {\bf pseudoramified algebraic dimension\/} of 
$M_{v_i}$ is equal to its rank $n_i=i ^n\cdot N^n$~.

Note that the geometric dimension of the $T_n(L_{v_i})$-subsemimodule $M_{v_i}$ is equal to ``~$n$~''.  So, the geometric and algebraic dimensions generally do not coincide.
\end{defi} \vskip 11pt

\begin{propo} \quad Let $\widetilde M^{(nr)}_R\otimes \widetilde M^{(nr)}_L$ denote the pseudo(un)ramified bisemisheaf over the real $GL_n(L^{(nr)}_{\o v}\times L^{(nr)}_{v})$-bisemimodule $(M^{(nr)}_R\otimes M^{(nr)}_L)$
isomorphic to its analytic counterpart $(M^{s(nr)}_R\otimes M^{s(nr)}_L)$~.

Let $(T^D_R\otimes T^D_L)$ be a differential bioperator acting on $(\widetilde M^{(nr)}_R\otimes \widetilde M^{(nr)}_L)$  and transforming them into the corresponding shifted bisemisheaves $(\widetilde M^{(nr)}_{R_{n[m]}}\otimes \widetilde M^{(nr)}_{L_{n[m]}})$~.

%Let $(\Pra\Gamma ^h_R\times \Pra \Gamma _L)$ be the bilinear subsemigroup of inner automorphisms of\linebreak $GL_m(\CC\times \CC)$ and let $(\Gamma ^h_R\times \Gamma ^h_L)$ be the bilinear semigroup of modular automorphisms of $GL_m(\CC\times \CC)$ such that $GL_m(\CC\times \CC)$ acts by cross product on $GL_n(\Aa_R\times \Aa_L)$~.

Then, the   bisemimodules $(M^{(nr)}_R\otimes M^{(nr)}_L)$ as 
well as their shifted  counterparts $(M^{(nr)}_{R_{n[m]}}\otimes M^{(nr)}_{L_{n[m]}})$  are characterized by the following ranks or algebraic dimensions:
\Bean
\item the pseudounramified bisemimodule $M^{nr}_R\otimes M^{nr}_L$ has for algebraic dimension $d=\sum\limits_{i=1}^d i^{n^2}$~;
\item pseudounramified shifted bisemimodule $M^{nr}_{R_{n[m]}}\otimes M^{nr}_{L_{n[m]}}$ has for algebraic dimension $d=\sum\limits_{i=1}^d i^{n^2}$ and for shifted algebraic dimension
$d_s=\sum\limits_{i=1}^q i^{m^2}$~;
\item pseudoramified bisemimodule $M_R\otimes M_L$ has for algebraic dimension $d=\sum\limits_{i=1}^d (i\cdot N)^{n^2}$~;
\item pseudoramified shifted bisemimodule  $M_{R_{[m]}}\otimes M_{L_{n[m]}}$ has for algebraic dimension $d=\sum\limits_{i=1}^d (i\cdot N)^{n^2}$ and for shifted algebraic dimension  $d_s=\sum\limits_{i=1}^q(i\cdot N)^{m^2}$~.
\Ee
\end{propo}
\vskip 11pt

\bpr This results from sections 3.5 and 3.6 and from \cite{30}.\epr
\vskip 11pt

%%%%%%%%%on a retiré ici § 3.7 et supprmé le paragraphe juste avant

\begin{propo} \quad  Under the ``action'' of the bioperator $(T^D_R\otimes T^D_L)\in 
\MM\RL(H^{\pm}_a)$ of rank $(m\times m)$~, the Shimura bisemivariety
\[\partial \o S_{G_{R\times L}} 
= P_n({L^+_{\o v^1}}\times {L^+_{v^1}})\setminus GL_n(L^+_{R_d}\times L^+_{L_d})\big/GL_n((\ZZ/N \ZZ)^2)\]
is shifted into $(m\times m)$ dimensions according to:
\[ T^D_R\otimes T^D_L : \quad \partial \o S_{G_{R\times L}}\To \partial \o S_{G_{R\times L;n[m]}}\]
where $\partial \o S_{G_{R\times L;n[m]}}$ is the shifted Shimura bisemivariety given by:
\begin{align*}
\partial \o S_{G_{R\times L;n[m]}}&=P_{n[m]}(({L_{\o v^1}}\otimes \rit)\times ( {L_{v^1}}\otimes \rit))\setminus\\
& \qquad \quad GL_{n[m]}((L^+_{R_d}\otimes \rit)\times (L^+_{L_d}\otimes \rit))\big/ GL_{n[m]}((\ZZ/N\ \ZZ)^2\otimes\rit^2)\;.\end{align*}
\end{propo} \vskip 11pt

\begin{propo} \quad The  bilinear cohomology semigroup of the Shimura bisemivariety 
$\partial \o S_{G_{R\times L}}$ is shifted under the action of the 
differential bioperator $(T^D_R\otimes T^D_L)\in\MM\RL(H^{\pm}_a)$ according to:
\begin{align*}
T^D_R\otimes T^D_L :\quad 
& H^{2j}(\partial\o S_{R\times L},\widetilde M^{2j}_R\otimes \widetilde M^{2j}_L) \\
& \quad \To 
H^{2j-2k}(\partial \o S_{G_{R\times L;n[m]}},\widetilde M^{}_{R_{2j[2k]}}\otimes \widetilde M^{2j}_{L_{2j[2k]}}) \end{align*}
in such a way that the shifted bilinear Eisenstein cohomology decomposes according to the bicosets of the quotient bisemigroupoid
\[ GL_{2j[2k]}((L_{\o v}\otimes \rit)\times (L_{v}\otimes \rit))\big/ GL_{2j[2k]}((\ZZ/N\ \ZZ)^2\otimes\rit^2)\]
as follows:
\[ H^{2j-2k}(\partial \o S_{G_{R\times L;n[m]}},\widetilde M^{2j}_{R_{2j[2k]_\oplus}}\otimes \widetilde M^{2j}_{L_{2j[2k]_\oplus}})=
\txt\bigoplus\limits^q_{i=1}\bigoplus\limits_{m_i} (\widetilde M^{2j}_{R_{2j[2k]}}
(i;m_i)\otimes \widetilde M^{2j}_{L_{2j[2k]}}(i;m_i))\]
where $m_i$ refers to the multiplicity of the shifted subbisemimodule 
$(M^{2j}_{R_{2j[2k]}}(i;m_i)\otimes\linebreak M^{2j}_{L_{2j[2k]}}(i;m_i))$~.
\end{propo} \vskip 11pt

\bpr According to proposition 2.3 and the Langlands bilinear global program developed in \cite{29} and in \cite{xxx}, we have that
\begin{align*}
H^{2j}(\partial \o S_{G_{R\times L}},\widetilde M^{2j}_{R_\oplus}\otimes \widetilde M^{2j}_{L_\oplus}) &\approx \FRepsp(GL_{2j}(L_{\o v_\oplus}\times L_{v_\oplus})\\
&= \txt\bigoplus\limits_{i=1}^q\bigoplus\limits_{m_i}(\widetilde M^{2j}_{\o v_{i;m_i}}\otimes \widetilde M^{2j}_{v_{i;m_i}})\;.
\end{align*}
Then, the shifted bilinear Eisenstein cohomology verifies:
\begin{align*}
H^{2j-2k}(\partial \o S_{G_{R\times L;n[m]}},\widetilde M^{2j}_{R_{2j[2k]_\oplus}}\otimes \widetilde M^{2j}_{L_{2j[2k]_\oplus}})
&\simeq \FRepsp (GL_{2j[2k]}((L_{\o v_\oplus}\otimes \rit)\times (L_{v_\oplus}\otimes \rit)))\\
&= \txt\bigoplus\limits^q_{i=1}\bigoplus\limits_{m_i}(\widetilde M^{2j}_{R_{2j[2k]}}(i;m_i)\otimes 
(\widetilde M^{2j}_{L_{2j[2k]}}(i;m_i))\end{align*}
such that:
\begin{equation} \FRepsp(GL_{2j[2k]}((L_{\o v_i}\otimes \rit)\times (L_{v_i}\otimes \rit))) = \widetilde M^{2j}_{R_{2j[2k]}}(i;m_i)\otimes \widetilde M^{2j}_{L_{2j[2k]}}(i;m_i))\;.
\tag*{\eop}\end{equation} \vskip 11pt

\begin{propo} \quad Let us fix the integers
\[ 1\le \ell_{R,L}\le j\;, \quad 1\le k_R\le j \quad \text{and} \quad 1\le u_{R,L}\le k\;, \quad 1\le v_R\le k\]
with the condition that $m\le n$~.

Then, the shifted bilinear Eisenstein cohomology decomposes into the direct sum of completely irreducible orthogonal or nonorthogonal shifted bilinear Eisenstein cohomologies according to:\\[6pt]
\textbullet\ $H^{2j-2k}(\partial \o S_{G\RL;n[m]},\widetilde M^{2j}_{R_{2j[2k]_\oplus}}\otimes_D\widetilde M^{2j}_{L_{2j[2k]_\oplus}})$
\\[6pt]
\mbox{}\quad $=\bigoplus\limits^q_{i=1}\bigoplus\limits_{m_i}\bigoplus\limits_{\ell_R=\ell_L}\bigoplus\limits_{u_{R,L}}
H^{2_{\ell_R}-2_{u_R}}(\partial \o S_{G\RL,n[m]},\widetilde M^2_{R_{2_{\ell_R}[2_{u_R}]}}(i;m_i)\otimes \widetilde M^2_{L_{2_{\ell_L}[2_{u_L}]}}(i;m_i))$
\\[6pt]
\mbox{}\quad $= \FRepsp(GL_{2j=2_1+\cdots+2_\ell+\cdots+2_{j[k]}}((L_{\o v}\otimes \rit)\times (L_{v}\otimes \rit))$ \vskip 11pt

\noindent \textbullet\ $H^{2j-2k}(\partial \o S_{G\RL;n[m]},\widetilde M^{2j}_{R_{2j[2k]_\oplus}}\otimes \widetilde M^{2j}_{L_{2j[2k]_\oplus}})$
\\[6pt]
\mbox{}\quad $=\bigoplus\limits^q_{i=1}\bigoplus\limits_{m_i}\bigoplus\limits_{\ell_R=\ell_L}\bigoplus\limits_{u_{R,L}}
H^{2_{\ell_R}-2_{u_R}}(\partial \o S _{G\RL;n[m]},\widetilde M^2_{R_{2_{\ell_R}[2_{u_R}]}}(i;m_i)\otimes \widetilde M^2_{L_{2_{\ell_L[2_{u_L}]}}}(i;m_i))$
\\[6pt]
\mbox{}\quad \quad $\bigoplus\limits^q_{i=1}\bigoplus\limits_{m_i}\bigoplus\limits_{k_R\neq\ell_L}\bigoplus\limits_{v_R\neq u_L}
H^{2_{k_R}-2_{v_R}}(\partial \o S _{G\RL;n[m]},\widetilde M^2_{R_{2_{k_R}[2_{v_R}]}}(i;m_i)\otimes \widetilde M^2_{L_{2_{\ell_L}[2_{u_L}]}}(i;m_i))$
\\[6pt]
\mbox{}\quad $= \FRepsp(GL_{2j\RL[2k]}((L_{\o v}\otimes \rit)\times (L_{v}\otimes \rit))$ \vskip 11pt

\noindent where $\widetilde M^2_{2_{\ell_L}[2_{u_L}]}(i;m_i))$ is a two-dimensional shifted functional representation space over the $T_{2_{\ell_L[u_L]}}(L_{v_i}\otimes \rit)$-semimodule.
\end{propo} \vskip 11pt

\bpr This proposition introduces the complete reducibility of the bilinear Eisenstein shifted cohomology semigroup in complete analogy with the unshifted case developed in proposition 2.6 and accordig to \cite{xxx}.\epr
 \vskip 11pt

%%%%%%%%%%coupe ici

\begin{defi} \quad {\bf Solvable bilinear Hilbert semispaces:\/}  \quad
\Be
\item Let $\widetilde M^{(nr)}_R\otimes \widetilde M^{(nr)}_L=\{\widetilde M^{(nr)}_{\o v_{i,m_i}}\otimes \widetilde M^{(nr)}_{v_{i,m_i}}\}^{q}_{i=1}$ be the bisemisheaf of differentiable bifunctions
\[\widetilde M^{(nr)}_{\o v_{i,m_i}}\otimes \widetilde M^{(nr)}_{v_{i,m_i}}\equiv 
\phi _{G_{i,m_{i_R}}}(x_{i_R})\otimes \phi _{G_{i,m_{i_L}}}(x_{i_L})\]
over the $GL_n(L^{(nr)}_{\o v}\times L^{(nr)}_v)$-bisemimodule $M^{(nr)}_R\otimes M^{(nr)}_L$ in such a way that
 $\widetilde M^{(nr)}_{R(P)}\otimes \widetilde M^{(nr)}_L$ is an extended internal (pseudounramified) bilinear Hilbert semispace $H^{+,(nr)}_a$ according to definitions 1.8.

Then, the $i$-th class $\{\widetilde M^{(nr)}_{\o v_{i,m_i}}\otimes \widetilde M^{(nr)}_{v_{i,m_i}}\}_{m_i}$ of 
$\widetilde M^{(nr)}_R \otimes \widetilde M^{(nr)}_L$ corresponds to the extended internal bilinear Hilbert subsemispace
$H^{+,(nr)}_a(i)$~; so that we get the towers
\begin{gather*}
H^{+,nr}_a(1) \subset \cdots \subset H^{+,nr}_a(i) \subset \cdots \subset H^{+,nr}_a(q)\;, \\
H^+_a(1) \subset \cdots \subset H^+_a(i) \subset \cdots \subset H^+_a(q)\;, \end{gather*}
of embedded pseudounramified and pseudoramified bilinear Hilbert subsemispaces.

Taking into account the isomorphism between the algebraic and analytic bilinear Hilbert semispaces $H^+_a$ and $H^+_h$~, corresponding  towers of embedded analytic bilinear Hilbert subsemispaces can also be envisaged:
\begin{gather*}
H^{+,nr}_h(1) \subset \cdots \subset H^{+,nr}_h(i) \subset \cdots \subset H^{+,nr}_h(q)\;, \\
H^+_h(1) \subset \cdots \subset H^+_h(i) \subset \cdots \subset H^+_h(q)\;. \end{gather*}

\item Let $\widetilde M^{nr}_{R_{L_\oplus}}\otimes_{(D)}\widetilde M^{nr}_{L_\oplus}=\bigoplus\limits^q_{i=1}(\widetilde M^{nr}_{\o v_i:m_i}\otimes_{(D)}\widetilde M^{nr}_{v_i;m_i})$ be the decomposition of the bisemisheaf over the  $GL_n(L^{nr}_{\o v_\oplus}\times L^{nr}_{v_\oplus})$-bisemimodule $M^{nr}_{R_{L_\oplus}}\otimes_{(D)}M^{nr}_{L_\oplus}$~.  Then, the algebraic pseudounramified extended (resp. diagonal) bilinear Hilbert semispace $H^{+,nr}_{a_\oplus}$ (resp. $\Hs^{+,nr}_{a_\oplus}$~) decomposes  according to:
\[ H^{+,nr}_{a_\oplus} = \txt\bigoplus\limits^q_{i=1}H^{+,nr}_a(i) \qquad \text{(resp.}\quad 
\Hs^{+,nr}_{a_\oplus}=\bigoplus\limits^q_{i=1}\Hs^{+,nr}_a(i)\ )\]
where $\widetilde M^{nr}_{\o v_i;m_i}\otimes \widetilde M^{nr}_{v_i;m_i}\simeq H^{+,nr}_a(i)$~.

So, we can construct {\bbf a tower of direct sums\/} of embedded algebraic pseudounramified extended (resp. diagonal) bilinear Hilbert subsemispaces
\[H^{+,nr}_a\{1\}\subset \cdots \subset H^{+,nr}_a\{i\}\subset \cdots\subset H^{+,nr}_a\{q\}\]
such that:
\Bi
\item $H^{+,nr}_a\{q\}\equiv H^{+,nr}_{a_\oplus}=\bigoplus\limits^q_{\nu =1}H^{+,nr}_a(\nu )$~,
\item \quad $H^{+,nr}_a\{i\} =\bigoplus\limits^i_{\nu=1}H^{+,nr}_a(\nu)$~,
\item \quad $\Hs^{+,nr}_a\{i\} =\bigoplus\limits^i_{\nu=1}\Hs^{+,nr}_a(\nu)$~,
\Ei
refer respectively to the $q$-th, $i$-th and $i$-th {\bf state\/} of $H^+_a$~, $H^{+,nr}_a$ and $\Hs^{,nr}+_a$~.

\item Considering the isomorphism between the algebraic and analytic bisemimodules\linebreak $(M^{nr}_R\otimes _{(D)} M^{nr}_L)$ and $(M^{s,nr}_R\otimes_{(D)}M^{s,nr}_L)$~, a tower of direct sums of embedded analytic pseudounramified extended (resp. diagonal) bilinear Hilbert subsemispaces also exists:
\begin{align*}
H^{+,nr}_h\{1\}&\subset \cdots \subset H^{+,nr}_h\{i\}\subset \cdots\subset H^{+,nr}_h\{q\}\\
\noalign{\qquad \qquad and}
\Hs^{+,nr}_h\{1\}&\subset \cdots \subset \Hs^{+,nr}_h\{i\}\subset \cdots\subset \Hs^{+,nr}_h\{q\}\end{align*}
such that:
\Bi
\item $H^{+,nr}_h\{q\}\equiv H^{+,nr}_{h_\oplus}=\bigoplus\limits^q_{\nu =1}H^{+,nr}_h(\nu )$~,
\item \quad $H^{+,nr}_h\{i\} =\bigoplus\limits^i_{\nu=1}H^{+,nr}_h(\nu)$~,
\item \quad $\Hs^{+,nr}_h\{i\} =\bigoplus\limits^i_{\nu=1}\Hs^{+,nr}_h(\nu)$~,
\Ei
refer respectively to the $q$-th, $i$-th and $i$-th {\bf state\/} of $H^{+,nr}_h$~, $H^{+,nr}_h$ and $\Hs^{+,nr}_h$~.

\item If the decomposition of the pseudoramified bisemisheaf $\widetilde M_{R_L}\otimes_{(D)}\widetilde M_L$ over the  $GL_n(L_{\o v}\times L_v)$-bisemimodule $M_{R_L}\otimes_{(D)}M_L$ is envisaged, then the algebraic  pseudoramified bilinear Hilbert semispaces $H^+_a$ and $\Hs^+_a$ decompose according to:
\Bi
\item \quad $H^+_a\{i\}=\bigoplus\limits^i_{j=1}H^+_a(j)$~;
\item \quad $\Hs^+_a\{i\}=\bigoplus\limits^i_{j=1}\Hs^+_a(j)$~;\quad $1\le j\le i$~,
\Ei
where $1\le i\le q$~.  This leads to {\bbf towers of direct sums\/} of embedded algebraic and analytic pseudoramified bilinear Hilbert subsemispaces, i.e. {\bf towers of states of these bilinear Hilbert semispaces\/} $H^+_a$~, $H^+_h$~, $\Hs^+_a$ and $\Hs^+_h$~:
\Bi
\item \quad $H^+_a\{1\} \subset \cdots \subset H^+_a\{i\} \subset \cdots \subset H^+_a\{q\}$~,
\item \quad $H^+_h\{1\} \subset \cdots \subset H^+_h\{i\} \subset \cdots \subset H^+_h\{q\}$~,
\item \quad $\Hs^+_a\{1\} \subset \cdots \subset \Hs^+_a\{i\} \subset \cdots \subset \Hs^+_a\{q\}$~,
\item \quad $\Hs^+_h\{1\} \subset \cdots \subset \Hs^+_h\{i\} \subset \cdots \subset \Hs^+_h\{q\}$~,
\Ei
where  $H^+_a\{1\} \equiv H^+_a(1) $~.

 The towers of embedded bilinear Hilbert subsemispaces lead to consider that these bilinear Hilbert semispaces are ``solvable'' and thus graded.
\Ee
\end{defi} \vskip 11pt

\begin{defi} \quad {\bf Projectors:\/} \quad {\bf (a)} \quad Let
\begin{gather*}
H^{+,nr}_a(1) \subset \cdots \subset H^{+,nr}_q(i) \subset \cdots \subset H^{+,nr}_a(q)\;, \\
H^+_a(1) \subset \cdots \subset H^+_a(i) \subset \cdots \subset H^+_a(q)\; \end{gather*}
be the two towers of embedded pseudounramified and pseudoramified bilinear Hilbert subsemispaces introduced in section 3.12.

Then, the following projectors:
\begin{align*}
P^{{\rm fac}(nr)}_{i\RL} : \quad H^{+,nr}_a(q) &\To H^{+,nr}_a(i)\;, \quad \forall\ 1\le i\le q\;, \\
P^{\rm fac}_{i\RL} : \quad H^+_a(q) &\To H^+_a(i)\;, \end{align*}
can be introduced, as it is done classically, in such a way that:
\Bi
\item $P^{{\rm fac}(nr)}_{i\RL} $ projects $ H^{+,nr}_a(q) $ onto the $i$-th pseudounramified bilinear Hilbert subsemispace $H^{+,nr}_a(i)$~;
\item $P^{\rm fac}_{i\RL} $ projects $ H^+_a(q) $ onto the $i$-th pseudoramified bilinear Hilbert subsemispace $H^+_a(i)$~.
\Ei
\vskip 11pt

\noindent{\bf (b)}\quad
Let $H^+_{a_\oplus}$ be an  extended bilinear Hilbert semispace decomposing according to:

\qquad \begin{minipage}{12cm}
\Bi
\item $H^+_{a_\oplus}=\bigoplus\limits^q_{i=1}H^+_a\{i\}$  such that $H^+_a\{i\}=\bigoplus\limits^i_{\nu=1}H^+_a(\nu)$~;
\item[or \textbullet]   $H^{+,nr}_{a_\oplus}=\bigoplus\limits^q_{i=1}H^{+,nr}_a\{i\}$  such that $H^{+,nr}_a\{i\}=\bigoplus\limits^i_{j=1}H^{+,nr}_a(j)$~.
\Ei\end{minipage}

Then, we can define the (bi)projectors of states:
\begin{align*}
P_{i_{R\times L}}^{nr}: \quad H^{+,nr}_{a_\oplus} &\To H^{+,nr}_a\{i\}\; , \tag*{$1\le i\le q$~,}\\
P_{i_{R\times L}}: \quad H^{+\phantom{,nr}}_{a_\oplus}  &\To H^+_a\{i\}\; , \end{align*}
mapping $H^{+,nr}_{a_\oplus}$ respectively into its closed extended bilinear subsemispace $H^{+,nr}_a\{i\}$ which is the $i$-th (bisemi)state.

The (bi)projectors $P_{i_{R\times L}}^{nr}$ and $P_{i_{R\times L}}$ are idempotent (bi)operators in such a way that the mappings they generate are  inverse   deformations (of Galois representations) as proved by the author elsewhere \cite{xx}.
\end{defi}

\begin{propo}\quad The operator $T^D_{R,L}(\Pra \Gamma ^{[m]}_{R,L})$ (resp. $T^D_{R,L}(\Gamma ^{[m]}_{R,L})$~) is a
random operator decomposing into a set of operators $\{T^D_{R,L}(g^{[m]}_{R,L}(i))\}_i$ (resp. $\{T^D_{R,L}(\gamma^{h[m]}_{R,L}(i))\}_i$~), $\forall\ 1\le i\le q$   according to the shifted pseudoramified (resp. pseudounramified) conjugacy classes of $GL_{n[m]}(L^{nr}_v\otimes\rit)$  associated with the $T^{(t)}_m(\rit)$-principal bundle.
\end{propo}
\vskip 11pt

\paragraph{Proof:}  Indeed, according to section 3.5,a random operator $T^D_{R,L}(\Pra \Gamma ^{[m]}_{R,L})=\linebreak
\{T^D_{R,L}(g^{[m]}_{R,L}(i))\}^q_{i=1}$ \quad (resp. $
T^D_{R,L}(\Gamma ^{[m]}_{R,L})=\{T^D_{R,L}(\gamma^{[m]}_{R,L}(i))\}^q_{i=1}$~), acting on an
extended bilinear Hilbert semispace $H^{\mp(nr)}_{a}$~, is a set
$\{T^D_{R,L}(g^{[m]}_{R,L}(i))\}^q_{i=1}\in\{\MM_{R,L}(H^{\mp}_{a}(i))\}$  \quad (resp. $\{T^D_{R,L}(\gamma ^{[m]}_{R,L}(i))\}^q_{i=1}\in \{\MM_{R,L}(H^{\mp,nr}_{a}(i))\}$~) such that the bilinear form:
\nobeqn t_R(\ell,m)&=& (T^D_R(g^{[m]}_R(\ell))\
e^\ell_{R},e^m_{L}) \\
\mbox{(resp.}\quad t_R(\mu,\nu)&=& (T^D_R(\gamma^{[m]}_R(\mu ))\
e^\mu _{R},e^\nu _{L})\;) \\
\noalign{or}
t_L(\ell,m)&=& (e^\ell_{R},T^D_L(g^{[m]}_L(m))\
e^m_{L})\\
\mbox{(resp.}\quad t_L(\mu ,\nu )&=& (e^\mu _{R},T^D_L(\gamma^{[m]}_L(\nu ))\
e^\nu _{L})\;)
\noeeqn
be measurable.
\vskip 11pt

$\MM_{R,L}(H^{\mp}_{a}(i))$ (resp. $\MM_{R,L}(H^{\mp,nr}_{a}(i))$~) is a von Neumann
subsemialgebra relative to  bounded operators on a closed connected
subsemispace $H^{\mp,(nr)}_{a}(i)$ of $H^{\pm,(nr)}_a$ referring to the $i$-th   conjugacy class of $GL_n(L^{(nr)}_{\o v}\times L^{(nr)}_v)$~.

These considerations are made in complete analogy with what is known for
random operators on linear Hilbert (semi)spaces \cite{9}.
\epr
\vskip 11pt

\begin{propo} \quad {\bf 1)\/}  Let $T^D_{R,L}(g^{[m]}_{R,L}(u))$ and $T^D_{R,L}(g^{[m]}_{R,L}(v))$ be two right or left  random operators such that $u<v$~.  Then, the  random operator $T^D_{R,L}(g^{[m]}_{R,L}(v))$ is an ``extension'' of the  random operator $T^D_{R,L}(g^{[m]}_{R,L}(u))$ corresponding to a difference of  conjugacy classes $(v-u)$~.

\noindent {\bf 2)\/} Let $T^D_{R,L}(\gamma ^{[m]}_{R,L}(o))$ and $T^D_{R,L}(\gamma ^{[m]}_{R,L}(p))$ be two right or left  random operators such that $o<p$~.  Then, $T^D_{R,L}(\gamma ^{[m]}_{R,L}(p))$ is an ``extension'' of $T^D_{R,L}(\gamma ^{[m]}_{R,L}(o))$ corresponding to a difference of  conjugacy classes $(p-o)$~.
\end{propo} \vskip 11pt

\begin{defi} \quad {\bf Towers of pseudoramified and pseudounramified von Neumann subsemialgebras:\/}  \quad {\bf (a)} \quad 
In connection with the definition 3.13 introducing towers of direct sums of embedded bilinear Hilbert subsemispaces, we shall define here towers 
of sums of  random operators:
\begin{align*}
T^D_{R,L}(g^{[m]}_{R,L}\{i\}) &= \txt\bigoplus\limits^i_{j=1}T^D_{R,L}(g^{[m]}_{R,L}(j)) \\
\text{(resp.} \quad 
T^D_{R,L}(\gamma ^{[m]}_{R,L}\{i\}) &= \txt\bigoplus\limits^i_{j =1}T^D_{R,L}(\gamma ^{[m]}_{R,L}(j ))\;), \end{align*}
such that 
\begin{alignat*}{3}
T^D_{R,L}(g^{[m]}_{R,L}\{i\})&\in \MM_{R,L}(H^{\mp}_a\{i\})\;, \quad &&1\le i\le q\;, \\
\text{(resp.}\quad
T^D_{R,L}(\gamma ^{[m]}_{R,L}\{i\})&\in \MM_{R,L}(H^{\mp,nr}_a\{i\})\;), \end{alignat*}
where $\MM_{R,L}(H^{\mp}_a\{i\})$ (resp. $\MM_{R,L}(H^{\mp,nr}_a\{i\})$~) is the pseudoramified (resp. pseudounramified) von Neumann subsemialgebra of the $i$-th state referring to the $i$-th sum of  random operators.

So, a tower of pseudoramified and pseudounramified von Neumann subsemialgebras of states can be introduced by:
\begin{align*}
\MM_{R,L}(H^{\mp}_a\{1\}) &\subset \cdots \subset  
\MM_{R,L}(H^{\mp}_a\{i\}) \subset \cdots \subset 
\MM_{R,L}(H^{\mp}_a\{q\})  \;, \\
\text{(resp.} \;
\MM_{R,L}(H^{\mp,nr}_a\{1\}) &\subset \cdots \subset  
\MM_{R,L}(H^{\mp,nr}_a\{i\}) \subset \cdots \subset 
\MM_{R,L}(H^{\mp,nr}_a\{q\})  \;), \end{align*}
such that
\begin{align*}
\MM_{R,L}(H^{\mp}_a\{i\}) &=\txt\bigoplus\limits^i_{j=1}
\MM_{R,L}(H^{\mp}_a(j)) \\
\text{(resp.} \quad
\MM_{R,L}(H^{\mp,nr}_a\{i\}) &=\txt\bigoplus\limits^i_{j =1}
\MM_{R,L}(H^{\mp,nr}_a(j ))\; ). \end{align*}
\vskip 11pt

\noindent {\bf (b)} \quad Similarly, on the towers
\begin{align*}
H^{\mp,nr}_h(1) & \subset \cdots \subset H^{\mp,nr}_h(i)  \subset \cdots \subset 
H^{\mp,nr}_h(q)  \\
\text{and} \quad
H^{\mp}_h(1) & \subset \cdots \subset H^{\mp}_h(i)  \subset \cdots \subset 
H^{\mp}_h(q) 
\end{align*}
of analytic Hilbert subsemispaces introduced in definition 3.12, the corresponding towers of pseudounramified and pseudoramified von Neumann subsemialgebras will be given by:
\begin{align*}
\MM_{R,L}(H^{\mp,nr}_h(1)) & 
\subset \cdots \subset \MM_{R,L}(H^{\mp,nr}_h(i))   \subset \cdots \subset 
\MM_{R,L}(H^{\mp,nr}_h(q))\;, 
\quad 1\le i\le q\le\infty \;, \\
\noalign{\noindent and by} 
\MM_{R,L}(H^{\mp}_h(1)) & 
\subset \cdots \subset \MM_{R,L}(H^{\mp}_h(i))   \subset \cdots \subset 
\MM_{R,L}(H^{\mp}_h(q))\;. 
\end{align*}
\end{defi} \vskip 11pt

\begin{propo}\quad Let $\MM_{R,L}(H^{\mp,(nr)}_{a})$ be the von
Neumann semialgebra of bounded self-adjoint operators on the smooth extended
bilinear Hilbert semispace $H^{\mp,(nr)}_{a}$~.

Let $\MM_{R,L}(H^{\mp,(nr)}_{a}\{i\})$ be the von Neumann
subsemialgebra of  random operators on the closed  smooth {\bf extended\/} bilinear
subsemispace  $H^{\mp,nr}_{a}\{i\}$ and let $\MM_{R,L}(\Hs^{\mp,nr}_{a}\{i\})$ be the corresponding von Neumann
subsemialgebra on the closed  smooth {\bf internal\/} diagonal bilinear subsemispace
$\Hs^{\mp,nr}_{a}\{i\}$~.

Then, the discrete spectrum $\sigma (T^D_{R,L})$ of an operator
$T^D_{R,L}\in \MM_{R,L}(H^{\mp,(nr)}_{a})$ is obtained by the
morphism:
\[ i^{a}_{\{i\}^D_{R,L}}\circ i^{a}_{\{i\}_{R,L}}\ : \quad
\begin{array}[t]{ccc}
\MM_{R,L}(H^{\mp,(nr)}_{a})&\To & [ \MM_{R,L}(\Hs^{\mp,(nr)}_{a}\{i\})]_i\;, \\
T^D_{R,L} &\To & \sigma (T^D_{R,L})\end{array}\]
where  $i^{a}_{\{i\}_{R,L}}$ and $i^{a}_{\{i\}^D_{R,L}}$ are given by:
\begin{alignat*}{3} 
i^{a}_{\{i\}_{R,L}}&:\quad \MM_{R,L}(H^{\mp,(nr)}_{a}) &&\To
[ \MM_{R,L}(H^{\mp,(nr)}_{a}\{i\})]_i\;, \\
\noalign{\vskip 6pt}
i^{a}_{\{i\}^D_{R,L}}&:\quad 
[ \MM_{R,L}(H^{\mp,(nr)}_{a}\{i\})]_i &&
\To [ \MM_{R,L}(\Hs^{\mp,(nr)}_{a}\{i\})]_i\;.\end{alignat*}
\end{propo}
\vskip 11pt

\paragraph{Proof:}  First remark that
$\MM_{R,L}(H^{\mp,(nr)}_{a})$ is a non-abelian von Neumann
semialgebra since the extended bilinear Hilbert semispace $H^{\mp,(nr)}_{a}$
constitutes the enveloping (semi)algebra of the semimodule $M^{(nr)}_L$ (resp.
$M^{(nr)}_{L_R}$~).

The morphism
\begin{align*} i^{a}_{\{i\}_{R,L}}\ : \quad 
 \MM_{R,L}(H^{\mp,(nr)}_{a}) &\To   [ \MM_{R,L}(H^{\mp,(nr)}_{a}\{i\})]_i\;, \\
T^D_{R,L}(\Pra\Gamma ^{[m]}_{R,L}) &\To   
[ T^D_{R,L}(g^{[m]}_{R,L}\{i\})]_i \\
 \mbox{(resp.} \quad
T^D_{R,L}(\Gamma ^{[m]}_{R,L})&\To   
[ T^D_{R,L}(\gamma^{[m]}_{R,L}\{i\}]_i\;),\end{align*}
transforms the bounded operator $T^D_{R,L}(\Pra\Gamma^{[m]}_{R,L})$ (resp. $T^D_{R,L}(\Gamma^{[m]}_{R,L})$~) into the
set\linebreak $[ T^D_{R,L}(g^{[m]}_{R,L}\{i\})]_i$
(resp. $[ T^D_{R,L}(\gamma^{[m]}_{R,L}\{i\}]_i$
  of bounded operators (i.e. sums of random operators acting on closed subsemispaces
$M_{L_{R,L}}\{i\}$ whose sums of enveloping subsemispaces are $H^{\mp,(nr)}_{h}\{i\}$~).

On the other hand, the isomorphism $i^{a}_{\{i\}^D_{R,L}}$
\[ i^{a}_{(i)^D_{R,L}}\ : \quad \begin{array}[t]{ccc}
[ \MM_{R,L}(H^{\mp,(nr)}_{a}\{i\})]_i&\To & 
[
 \MM_{R,L}(\Hs^{\mp,(nr)}_{a}\{i\}]_i\;, \\
{[T^D_{R,L}(g^{[m]}_{R,L}\{i\}]_i}
&\To 
&\sigma 
(T^D_{R,L})\;,\\
\mbox{(resp.\ }\quad 
[T^D_{R,L}(\gamma^{[m]}_{R,L}\{i\}]_i&\To &\sigma_{nr}
(T^D_{R,L})\;),\end{array}\]
transforms the non-abelian von Neumann subsemialgebra $ \MM_{R,L}(H^{\mp,(nr)}_{a}\{i\})$ into the abelian or diagonal von
Neumann subsemialgebra $\MM_{R,L}(\Hs^{\mp,(nr)}_{a}\{i\})$ of  sums of
random operators.  $[\MM_{R,L}(\Hs^{\mp,(nr)}_{a}\{i\}]_i$ is
then an algebra of the sum of random operators acting on diagonal enveloping
subsemialgebras $(\Hs^{\mp,(nr)}_{a}\{i\})$~.  $\sigma (T^D_{R,L})$ (resp. $\sigma _{nr}(T^D_{R,L})$~) is thus the pseudoramified (resp. pseudounramified) spectrum of
the bounded operator $T^D_{R,L}$~.\epr
\vskip 11pt

\begin{coro} Let $\MM_{R\times L}(H^{\mp,(nr)}_{a})$ be the
von Neumann bisemialgebra of bounded bioperators $T^D_R\otimes T^D_L$ on
$H^{\mp,(nr)}_{a}$ and let $\MM_{R\times L}(H^{\mp,(nr)}_{a}\{i\})$ be the $i$-th corresponding von Neumann subbisemialgebra of the sum of random bounded
bioperators on $H^{\mp,(nr)}_{a}\{i\}$~.  If $\MM^{h}_{R\times L}(\Hs^{\mp}_{a}\{i\})$
is the $i$-th von Neumann diagonal subbisemialgebra of random diagonal bioperators
$T^D_R\{i\}\otimes_DT^D_L\{i\}$ on $\Hs^{\mp,nr}_{a}\{i\}$~, then the
discrete spectrum $\sigma (T^D_R\otimes T^D_L)$ of $(T^D_R\otimes
T^D_L)\in \MM_{R\times L}(H^{\mp,(nr)}_{a})$ is obtained by
the morphism:
\[ i^{a}_{\{i\}^D_{R\times L}}
\circ i^{a}_{\{i\}_{R\times L}}
\ : \quad \begin{array}[t]{ccc}
\MM_{R\times L}(H^{\mp,(nr)}_{a})&\To & 
[ \MM_{R\times L}(\Hs^{\mp,(nr)}_{a}\{i\})]_i\;, \\
T^D_{R\times L}&\To &\sigma (
T^D_{R\times L})\;,\end{array}\]
where $T^D_{R\times L}$ is the condensed notation for $T^D_R\otimes
T^D_L$~.

$\MM_{R,L}(H^{\pm,(nr)}_a)$  then corresponds to a solvable (bi)semialgebra.
\end{coro}
\vskip 11pt

\paragraph{Proof:}  \quad This corollary is an extension of the preceding
proposition to the bioperator $(T^D_R(\Gamma _R)\otimes T^D_L(\Gamma _L))$~.\epr
\vskip 11pt

\sub{Shifted global pseudounramified (resp. pseudoramifed) elliptic bisemimodules}  Referring to proposition 3.6, 
the action of the differential bioperator $(T^D_R\otimes T^D_L)$ on the 
bisemisheaf $ ( \widetilde M^{(nr)}_R\otimes \widetilde M^{(nr)}_L)$ over the 
$GL_n(L^{(nr)}_{\o v}\times L^{(nr)}_v)$-bisemimodule $(M^{(nr)}_R\otimes M^{(nr)}_L)$ consists in mapping it into the shifted 
bisemisheaf $(\widetilde M^{(nr)}_{R_{n[m]}}\otimes \widetilde M^{(nr)}_{L_{n[m]}})$ over the
$GL_{n[m]}((L^{(nr)}_{\o v}\otimes \rit)\times (L^{(nr)}_v\otimes \rit))$-bisemimodule  
$(M^{(nr)}_{R_{n[m]}}\otimes M^{(nr)}_{L_{n[m]}})$ such that $(\widetilde M^{(nr)}_{R_{n[m]_\oplus}}\otimes \widetilde M^{(nr)}_{L_{n[m]_\oplus}})$ decomposes into   ``~$q$~'' subbisemisheaves.

But, according to proposition 2.10 referring to the Langlands global program introduced in \cite{29}, there is a bijection between the  $GL_n(L^{(nr)}_{\o v_\oplus}\times L^{(nr)}_{v_\oplus})$-bisemimodule\linebreak $(M^{(nr)}_{R_\oplus}\otimes M^{(nr)}_{L_\oplus})$ and its cuspidal counterpart given by the global pseudoramified (resp. pseudounramified) elliptic $G_{s\RL}$-bisemimodule:
\begin{align*}
\ELLIP\RL(n,q) &= \txt\sum\limits^q_{i=1}\sum\limits_{m_i}\lambda
(n,i,m_i)e^{-2\pi i(i)z} \otimes 
\sum\limits^q_{i=1}\sum\limits_{m_i}\lambda (n,i,m_i)e^{2\pi i(i)z}\;, 
\tag*{$z\in\RR^n$~,}\\
\text{(resp.} \quad
\ELLIP^{nr}\RL(n,q) &= \txt\sum\limits^q_{i=1}\sum\limits_{m_i}
\lambda _{nr}(n,i,m_i)e^{-2\pi i(i)z} \otimes 
\sum\limits^q_{i=1}\sum\limits_{m_i}\lambda _{nr}(n,i,m_i)e^{2\pi i(i)z}\;),
\end{align*}
such that we have the commutative diagram:
\[ \begin{CD}
\widetilde M^{nr}_R\otimes \widetilde M^{nr}_L @>{T^D_R\otimes T^D_L}>> \widetilde M^{nr}_{R_{n[m]}}\otimes \widetilde M^{nr}_{L_{n[m]}}\\
@VVV @VVV\\
\ELLIP^{nr}\RL(n,q) @>{T^D_R\otimes T^D_L}>> \ELLIP^{nr}\RL(n[m],q) \\
@VVV @VVV\\
\ELLIP\RL(n,q) @>{T^D_R\otimes T^D_L}>> \ELLIP\RL(n[m],q)\\
@AAA @AAA\\
\widetilde M_R\otimes \widetilde M_L @>{T^D_R\otimes T^D_L}>> \widetilde M_{R_{n[m]}}\otimes M_{L_{n[m]}}
 \end{CD}\]
where $\ELLIP^{nr}\RL(n[m],q)$ (resp. $\ELLIP^{nr}\RL(n[m],q)$~) is the shifted global pseudounramified (resp. pseudoramified) elliptic $((G_{s_R}\otimes \rit)\times (G_{s_L}\otimes \rit)$-bisemimodule. \vskip 11pt

As an application of proposition 3.17, we suggest the following proposition \cite{xxx}.
\vskip 11pt

\begin{propo} \quad The shifted global pseudounramified (resp. pseudoramified) $n$-\linebreak dimensional elliptic bisemimodule
\begin{align*}
\ELLIP^{nr}\RL(n[m],q)&= \ELLIP^{nr}_{R}(n[m],q)\otimes
\ELLIP^{nr}_{L}(n[m],q)\\
\text{(resp.} \quad
\ELLIP\RL(n[m],q)&= \ELLIP_{R}(n[m],q)\otimes
\ELLIP_{L}(n[m],q)\;),\end{align*}
gives rise to (or is functorially equivalent to) the eigenbivalue equation of the $i$-th (bi)states:
\begin{align*}
& (T^D_R\otimes T^D_L)( \ELLIP^{nr}_{R}(n,i)\otimes
\ELLIP^{nr}_{L}(n,i) \\
&\qquad = E^{nr}_R\{n,i\}\times E^{nr}_L\{n,i\}\cdot
 (\ELLIP^{nr}_{R}(n,i)\otimes
\ELLIP^{nr}_{L}(n,i)\;, && 1\le i\le q\;, \\
\text{(resp.} \quad
& (T^D_R\otimes T^D_L)( \ELLIP\RL(n,i)) \\
&\qquad = E_R\{n,i\}\times E_L(\{n,i\}\cdot
 (\ELLIP\RL(n,i)\;)\;, && 1\le i\le q\;.\end{align*}
\end{propo} \vskip 11pt

\bpr \Be \item The shifted global pseudounramified elliptic bisemimodule $\ELLIP^{(nr)}\RL(n[m],q)$ generates the eigenbivalue equation:
\[ \ELLIP^{nr}\RL(n[m],i)=(E^{nr}\RL\{n,i\})(\ELLIP^{nr}\RL(n,i))\]
which can be rewritten according to \cite{xxx}:
\[ (T^D_R\otimes T^D_L)(\ELLIP^{nr}\RL(n,i))=(E^{nr}_R\{n,i\}\times E^{nr}_L\{n,i)\}(\ELLIP^{nr}\RL(n,i))\]
where the \rl eigenvalue $E^{nr}_R(\{n,i\}$ (resp. $E^{nr}_L\{n,i\}$~) corresponds to a sum 
over the $i$ first pseudounramified algebraic classes
of shifts into $m$ dimensions of the Hecke characters $\lambda _{nr}(n,\nu ,m_\nu )$ (resp. $\lambda _{nr}(n,\nu ,m_\nu )$~) i.e. to infinitesimal generators of the considered Lie algebra, $1\le \nu \le i$~.

\item The bisemialgebra of von Neumann $\MM\RL(H^{\mp,nr}_h)$ can then be considered as a solvable bisemialgebra generating a tower of sums of pseudounramified von Neumann subbisemialgebras according to definition 3.16.  On the other hand, the set of pseudounramified eigenbivalues of $(T^D_R\otimes T^D_L)$ forms an embedded sequence:
\[ E^{nr}_R\{n,1\}\cdot E^{nr}_L\{n,1\}\subset \cdots 
E^{nr}_R\{n,i\}\cdot E^{nr}_L\{n,i\}\subset \cdots 
E^{nr}_R\{n,q\}\cdot E^{nr}_L\{n,q\}\]
in one-to-one correspondence with the set of embedded eigenbifunctions given by  the product, right by left, %of the $i$-th term 
of the truncated Fourier series at ``~$i$~'' terms:
\[ \ELLIP^{nr}\RL(n,i) = \txt\sum\limits^i_{\nu=1}\sum\limits_{m_\nu}\lambda _{nr}(n,\nu,m_\nu) e^{-2\pi i\nu z}\otimes \sum\limits^i_{\nu=1}\sum\limits_{m_\nu}\lambda _{nr}(n,\nu,m_\nu) e^{2\pi i\nu z}\ , \; z\in\RR^n\;. \]

\item The proof was given for the ``pseudounramified'' case, taking into account that the ``pseudoramified'' case can be handled similarly.\epr
\Ee

\begin{propo}\quad The discrete spectrum $ \sigma (T^D_{{R,L}\atop{\RL}})$ of 
$(T^D_{{R,L}\atop{\RL}}) \in \MM_{{R,L}\atop{\RL}}(H^{\mp,(nr)}_h)$ and the discrete spectrum
$ \sigma^a (T^D_{{R,L}\atop{\RL}})$ of $(T^D_{{R,L}\atop{\RL}}) \in \MM_{{R,L}\atop{\RL}}(H^{\mp,(nr)}_a)$ are isomorphic (and often equal).
\end{propo}

\paragraph{Proof:}  Consider the commutative diagram:
\[ \begin{CD}
\MM_{{R,L}\atop{R\times L}}(H^{\mp,(nr)}_{h})  @<{i_{\MM^a_{{R,L}\atop{R\times
L}}-\MM^{h}_{{R,L}\atop{R\times L}}}}<< \MM_{{R,L}\atop{R\times L}}(H^{\mp,(nr)}_{a}) \\
@V{i^{h}_{\{i\}_{{R,L}\atop{R\times L}}}}VV @VV{i^{a}_{\{i\}_{{R,L}\atop{R\times L}}}}V\\
[\MM_{{R,L}\atop{R\times L}}(H^{\mp,(nr)}_{h}\{i\})]_i @<{\quad i^{a-{h}}_{{R,L}\atop{R\times L}}\quad}<<
[\MM_{{R,L}\atop{R\times L}}(H^{\mp,(nr)}_a\{i\})]_i\\
@V{i^{h}_{\{i\}^D_{{R,L}\atop{R\times L}}}}VV @VV{i^{a}_{\{i\}^D_{{R,L}\atop{R\times L}}}}V\\
[\MM_{{R,L}\atop{R\times
L}}(\Hs^{\mp,(nr)}_{h}\{i\})]_i  @>{\quad i^{a-{h}}_{D_{{R,L}\atop{R\times L}}}\quad}>>
[\MM_{{R,L}\atop{R\times
L}}(\Hs^{\mp,(nr)}_a\{i\})]_i
\end{CD}\]
where
\Bi\item the isomorphism $i_{\MM^a_{{R,L}\atop{R\times L}}-\MM^{h}_{{R,L}\atop{R\times L}}}$ 
has been introduced in proposition 3.3;
\item the morphisms ${i^{h}_{\{i\}_{{R,L}\atop{R\times L}}}}$ and ${i^{a}_{\{i\} _{{R,L}\atop{R\times L}}}}$ result from the decomposition of $\MM_{{R,L}\atop{R\times L}}(H^{\mp,(nr)}_h)$ and of $\MM_{{R,L}\atop{R\times L}}(H^{\mp,(nr)}_a)$ into sums of pseudounramified or pseudoramified subbisemialgebras (see definition 3.16).
\Ei

From the isomorphism $i^{a-h}_{D_{{R,L}\atop{\RL}}}$~, it results that the
discrete spectrum $\sigma (T^D_{{R,L}\atop{R\times L}})$ of
$T^D_{{R,L}\atop{R\times L}}\in \MM_{{R,L}\atop{R\times
L}}(H^{\mp,(nr)}_{h})$ and the discrete spectrum $\sigma
^a(T^D_{{R,L}\atop{R\times L}})$ of $T^D_{{R,L}\atop{R\times
L}}\in \MM_{{R,L}\atop{R\times L}}(H^{\mp,(nr)}_a)$ are isomorphic. So, we get the thesis.\epr
\vskip 11pt

\sub{Factors of von Neumann}
\Bi
\item We are now interested in the classification of the factors of von Neumann, i.e. in von Neumann algebras having trivial centers (reduced to $\CC$~).  According to definition 3.16, we see that two types of towers of von Neumann sub(bi)semialgebras have been introduced:
\Bi
\item the first referring to  {\bf pseudounramified (algebraic) classes\/} of the bilinear Hilbert semispaces $H^{\mp,nr}_a$ (or $H^{\mp,nr}_h$~) on which they have been defined;
\item the second referring to  {\bf pseudoramified (algebraic) classes\/} of the bilinear Hilbert semispaces $H^{\mp}_a$ (or $H^{\mp}_h$~).
\Ei

So, the classification of factors of von Neumann will be based on these two types of towers of von Neumann subsemialgebras on bilinear  Hilbert (sub)semispaces which are associated with Hecke sublattices as developed in proposition 2.3 (proof).  As a result, the dimensions of the factors of von Neumann will directly refer to Hecke sublattices.

\item The bilinear Hilbert semispaces $H^{\mp}_a$~, isomorphic to $H^{\mp}_h$~, constituting the 
natural representation spaces of the von Neumann (bi)semialgebras, were supposed to be pseudoramified in the 
sense that the $GL_n(L_{\o v}\times L_v)$-bisemimodule $(M_{R_L}\otimes M_L)$ is pseudoramified.  That is to 
say that the  $T_n(L_{v_i})$-subsemimodule $M_{v_i}$ (as well as $M_{\o v_i}$~) has a rank given 
by $n_i=i^n\cdot N^n$ (see section 1.5).  

On the other hand, the corresponding pseudounramified 
$T_n(L^{nr}_{v_i})$-subsemimodule $M^{nr}_{v_i}$ would have a rank $n_i^{nr}=i^n$ according to 
\cite{29}, which allows to envisage the introduction of  pseudounramified bilinear Hilbert subsemispaces, noted $H^{nr}_a(i)$~, as it was defined in section 3.12.
\Ei \vskip 11pt

\begin{propo} {\bf (Classification of (bi)factors of von Neumann with respect to algebraic dimensions)}
\Be\item {\bbf Type ${\rm I}_i$\/}~: on the pseudounramified bilinear Hilbert semispace $H^{nr}_a$~, there are $q$ factors $\MM_{R,L}(H^{nr}_a(i))$ of type ${\rm I}_i$~, $1\le i\le q\le \infty$~, where $i$ denotes a global residue degree.

\item {\bbf Type ${\rm II}_1$\/}~: on the bilinear Hilbert subsemispace $H^{\mp,{\rm in}}_a[{L_{\o v^1}}\times {L_{v^1}}]$ restricted to the 
representation space of the bilinear parabolic subsemigroup $P_n({L_{\o v^1}}\times {L_{v^1}})$~, there are 
$N$ subfactors\linebreak $\MM_{R,L}(H^{\mp{\rm in}}_a(i))$~, $1\le i\le N$~, of type ${\rm I}{\rm I}_{1_i}$~, where $i$ denotes an internal algebraic dimension corresponding to the number of automorphisms of the global inertia subgroup.

The factor $\MM_{R,L}(H^{\mp,{\rm in}}_a(N))$ is the factor of type ${\rm I}{\rm I}_1$~.

\item {\bbf Type ${\rm II}_{(\infty)}$\/}~: on the tensor products $H^{\mp}_a(i)=H^{nr}_a(i)\otimes H^{\mp,{\rm in}}_a(N)$ of the 
pseudounramified bilinear Hilbert semispace $H^{nr}_a(i)$ by the bilinear Hilbert subsemispace $H^{\mp,{\rm in}}_a(N)$~, there are $q$ pseudoramified factors $\MM_{R,L}(H^{nr}_a(i)\otimes H^{\mp}_a(N))$ of type ${\rm II}_{(\infty)}$~, $1\le i\le q\le \infty$~, where $i$ denotes a global residue degree.

\item {\bbf Type ${\rm II}_\infty$\/}~: on the tensor products $H^{nr}_a(\infty )\otimes H^{\mp,{\rm in}}_a(j)$~, $1\le j\le N$~, the factors
$\MM_{R,L}(H^{nr}_a(\infty ))\otimes \MM_{R,L}(H^{\mp,{\rm in}}_a(j ))$~,  of type ${\rm I}{\rm I}_\infty $ are defined.
\Ee
\end{propo} \vskip 11pt

\bpr \nopagebreak \Be\item As there are $q$  conjugacy classes of the pseudounramified bilinear Hilbert semispace $H^{nr}_a$~, there are $q$ ``pseudounramified'' factors $\MM_{R,L}(H^{nr}_a(i))$ in the tower: 
\[ \MM_{R,L}(H^{nr}_a(1)) \subset \cdots \subset 
\MM_{R,L}(H^{nr}_a(i)) \subset \cdots \subset 
\MM_{R,L}(H^{nr}_a(q)) \]
as introduced in sections 3.12, 3.13 and  3.16. 

So, there are $q$ factors of type ${\rm I}_i$~, $1\le i\le q\le \infty$ with minimal projections.

\item If we consider the $N$ internal conjugacy classes of the bilinear parabolic semigroup $P_n(L_{\o v^1}\times L_{v^1})$ corresponding to the (shifted) intermediate inner automorphisms of the global inertia subgroups ${\rm I}_{L_{v_i}}$ having an order $N$~, we can introduce on $H^{\mp}_a[L_{\o v^1}\times L_{v^1}]$ a tower of inner hyperfinite subfactors \cite{23}, \cite{24}:
\[ \MM_{R,L}(H^{\mp,{\rm in}}_a(1)) \subset \cdots \subset 
\MM_{R,L}(H^{\mp,{\rm in}}_a(i)) \subset \cdots \subset 
\MM_{R,L}(H^{\mp,{\rm in}}_a(N)) \]
in such a way that:
\Bi
\item the index $[\MM_{R,L}(H^{\mp,{\rm in}}_a(i)):\MM_{R,L}(H^{\mp,{\rm in}}_a(1))]=i$ of the $i$-th hyperfinite subfactor with respect to the first hyperfinite subfactor is the internal algebraic dimension (see section 3.7).

\item the upper hyperfinite subfactor $\MM_{R,L}(H^{\mp,{\rm in}}_a(N))$ is the hyperfinite factor of type ${\rm I}{\rm I}_1$ having an index $N$ and corresponding to the order of the global inertia subgroup ${\rm I}_{L_{v_i}}$~.
\Ei

Indeed, if we take into account proposition 2.9, the Hecke characters on sublattices associated with $\MM_{R,L}(H^{\mp,{\rm in}}_a(i))$~, $1\le i\le N$~, must take values in the interval $[0,1]$~: they then correspond to the continued dimensions \cite{17} of the classes of the projectors of the subfactors of type ${\rm I}{\rm I}_1$ of von Neumann algebras on a linear Hilbert semispace.

\item As on the pseudounramified bilinear Hilbert subsemispaces $H^{nr}_a(i)$~, pseudo-\linebreak unramified factors 
$\MM_{R,L}(H^{nr}_a(i))$ of type ${\rm I}_i$ are defined and as, on the bilinear Hilbert subsemispace $H^{\mp,{\rm in}}_a(N)$~, a factor of type ${\rm I}{\rm I}_1$ is defined, it is evident that, on their tensor products $H^{nr}_a(i)\otimes H^{\mp,{\rm in}}_a(N)$~, pseudoramified factors of type ${\rm II}_i$~, characterizd by minimal projections, $1\le i\le q \le \infty$~, can be defined, the factor of type ${\rm I}{\rm I}_1$ ``ramifying'' the pseudounramified factors ${\rm I}_i$~.

\item And, then, the classical factors of Araki-Woods \cite{1}, \cite{15} of type ${\rm II}_{\infty}$ correspond to the factors $\MM_{R,L}(H^{nr}_a(i=\infty)) \otimes \MM_{R,L}(H^{\mp,{\rm in}}_a(j))$~, $1\le j\le N$~,  where
\Bi
\item $\MM_{R,L}(H^{nr}_a(i=\infty))$ is the pseudounramified factor of type ${\rm I}_{\infty}$~;
\item $\MM_{R,L}(H^{\mp,{\rm in}}_a(j))$ is the hyperfinite subfactor of type ${\rm II}_{1_j} $~.\epr
\Ei\Ee \vskip 11pt

\begin{coro} \quad The equivalent of a factor of type ${\rm I}{\rm I}{\rm I}_\lambda$ can be obtained by considering the cross product of the factor ${\rm I}{\rm I}_\infty$ by a subgroup of automorphisms of it \cite{1}, \cite{15}.
\end{coro} \vskip 11pt

\bpr  Indeed,  a factor $M_\lambda$ of type $M_\lambda$ \cite{28} is isomorphic to the cross product of a factor ``~$N$~'' of type ${\rm I}{\rm I}_\infty$ by $\Aut ``~N$~'' \cite{14}, \cite{16}, \cite{31}.\epr

\end{document}

%% file: cmds.tex
\def\simarrow{\mathrel{\raise -0.5mm\hbox{$\sim$}}\hspace{-1.8mm}{\rightarrow} } 

\def\bsimarrow{\leftarrow\hspace{-0.7mm}\mathrel{\raise -0.5mm\hbox{$\backsim$}} }

\def\bt{\begin{tabular}}
\def\te{\end{tabular}}

\def\lettrine#1#2#3{\noindent\hangindent#1\hangafter-#2
\hskip-#1\smash{\hbox to #1{#3\hfill}}\ignorespaces}

\newcommand{\To}[1]{\mathop{\to}\limits_{#1}}

\def\BM{\begin{pmatrix}}
\def\EM{\end{pmatrix}}

\def\txt{\textstyle}

\def\d=f{\buildrel\hbox{\scriptsize d\'{e}f}\over \Longleftrightarrow}

\def\cit{\text{\it I\hskip -6ptC\/}}

\def\square{\hfill\hbox{\vrule height .9ex width .8ex depth -.1ex}}
\def\rit{\text{\it I\hskip -2pt  R}}

\def\rl {\rit^{\hskip 1pt\ell}}

\def\Bd{{\text B}}

\def\Ed{{\text E}}

\def\Hs{{\cal H}}

\def\Ks{{\cal K}}

\def\Ls{{\cal L}}
\def\Ms{{\cal M}}

\def\be{\begin{equation}}
\def\ee{\end{equation}}
\def\beqn{\begin{eqnarray}}
\def\eeqn{\end{eqnarray}}
\def\nobeqn{\begin{eqnarray*}}
\def\noeeqn{\end{eqnarray*}}
\def\ba{\left(\begin{array}}
\def\ea{\end{array} \right) }

\def\bpr{\paragraph{Proof.}}

\def\epr{\square\vskip 6pt}
\def\eop{\hbox{\vrule height .9ex width .8ex depth -.1ex}}

\def\o{\overline}
\def\and{\; \mbox{and} \;}

\newcommand{\half}{\frac{1}{2}}

\def\hfl#1#2{\smash{\mathop{\hbox to 12mm{\rightarrowfill}}
\limits^{\scriptstyle #1}_{\scriptstyle #2}}}

\def\mod{\mathop{\rm mod}\nolimits}

\def\Be{\begin{enumerate}}
\def\Ee{\end{enumerate}}

\def\Bena{\begin{enumerate}
\def\labelenumi{\theenumi)}
\def\theenumi{\arabic{enumi}}
\def\labelenumii{\theenumii)}
\def\theenumii{\alph{enumii}}}

\def\Bean{\begin{enumerate}
\def\labelenumii{\theenumii)}
\def\theenumii{\arabic{enumii}}
\def\labelenumi{\theenumi)}
\def\theenumi{\alph{enumi}}}

\def\Bero{\begin{enumerate}
\def\labelenumii{\theenumii)}
\def\theenumii{\arabic{enumii}}
\def\labelenumi{(\theenumi)}
\def\theenumi{\roman{enumi}}}

\def\BeRo{\begin{enumerate}
\def\labelenumii{\theenumii)}
\def\theenumii{\arabic{enumii}}
\def\labelenumi{(\theenumi)}
\def\theenumi{\Roman{enumi}}}

\def\Bi{\vskip 11pt\begin{itemize}\itemsep=18pt}
\def\bi{\begin{itemize}}%\topsep=18pt\itemsep=18pt\partopsep=18pt\parsep=18pt}
\def\Ei{\end{itemize}\vskip 11pt}
\def\ei{\end{itemize}}

\def\Bd{\begin{description}}
\def\Ed{\end{description}}

\def\R{\right}
\def\L{\left}
\def\F{\frac}

%% file: Algebraic.bbl
\vskip 11pt



\vskip 22pt

\begin{thebibliography}{99}

\bibitem{1} {\sc Araki, H., Woods, E.J.:}  A classification of factors. {
Publ. Res. Inst. Math. Sc.\/} (Kyoto), {\bf 4\/} (1968), 51--130.

\bibitem{2} {\sc Artin, M.:}  On Azumaya algebras and finite dimensional
representations of rings. {\em J. of Algebras\/}, {\bf 11\/} (1969),
531--563.

\bibitem{3} {\sc Atiyah, M.F., Singer, I.:}  The index of elliptic operators
on compact manifolds.  {\em Bull. Amer. Math. Soc.\/}, {\bf 69\/} (1963),
422--433.

\bibitem{4}  {\sc Atiyah, M.F.:}  $K$-theory.  Benjamin,  1967.

\bibitem{5} {\sc Baum, P., Fulton, W., Mac Pherson, R.:}  Riemann-Roch
and topological $K$-theory for singular varieties.  {\em Acta Math.\/}, {\bf
143\/} (1979), 155--192.

\bibitem{6} {\sc Borel, A.:}  Linear algebraic groups.  {\em Grad. Texts in
Math.\/}, {\bf 126\/} (1991), Springer Verlag.


\bibitem{7} {\sc Borel, A.:}  Groupes linéaire algébriques.  {\em Annals of
Math.\/}, {\bf 64\/} (1956), 20--82.

\bibitem{8} {\sc Borel, A.:}  Regularization theorems in Lie algebra
cohomology.  {\em Duke Math. J.\/}, {\bf 50\/} (1983), 605--623.


\bibitem{9} {\sc Bratteli, O., Robinson, D.:}  Operator algebras and quantum statistical mechanics.  1.  (1979), Springer.

\bibitem{10} {\sc Brown, L.:}  Operator algebras and algebraic $K$-theory. 
{\em Bull. Amer. Math. Soc.\/}, {\bf 81\/} (1975), 1119--1121.

\bibitem{11} {\sc Brown, L., Douglas R., Fillmore, P.:} 
Extensions of $C^*$-algebras and $K$-homology.  {\em Annals of Math.\/},
{\bf 105\/} (1977), 265--324.

\bibitem{12} {\sc Casselman, W.:} Introduction to the Schwartz space of $\Gamma 
\setminus G$~.  {\em Can. J. Math.\/}, {\bf 15\/} (1989), 285--320.

\bibitem{13} {\sc Connell, J.C., Robson, J.C.:}  Non commutative Noetherian
rings.  Wiley, 1987.

\bibitem{14} {\sc Connes, A.:} Une classification des facteurs de type III. 
{\em Ann. Scient. \'Ecole Norm. Sup.\/}, {\bf 6\/} (1973), 133--252.

\bibitem{15} {\sc Connes, A.:} Sur la classification des facteurs de type
II.  {\em C.R. Acad. Sci.\/} (Paris), {\bf A281\/} (1975), 13--15.

\bibitem{16} {\sc Connes, A.:}  Classification of injective factors, Cases
$II_1$~, $II_\infty $~, $III_\lambda $~.  {\em Annals of Math.\/} {\bf
104\/} (1976), 73--115.

\bibitem{17} {\sc Elliott, G.A.:}  On the classification of inductive limits
of sequence of semisimple finite dimensional algebras.  {\em J. of Alg.\/}
{\bf 38\/} (1976), 29--44.

\bibitem{18} {\sc Farb, B., Dennis, K.:}  Non commutative algebras. {\em
Grad. Texts in Math.\/}, {\bf 144\/} (1993), Springer Verlag.

\bibitem{19} {\sc Feldman, J., Moore, C.:}  Ergodic equivalence relations,
cohomology and von Neumann algebras.  {\em Trans. Amer. Math. Soc.\/},
{\bf 234\/} (1974), 289--324, 325--359.

\bibitem{20} {\sc Grothendieck, A.:}  On the de Rham cohomology of algebraic
varieties.  {\em Publ. Math. IHES\/}, {\bf 29\/} (1966), 351--359.

\bibitem{21} {\sc Harder, G.:} Eisenstein cohomology of arithmetic groups. The
case GL2.  {\em Invent. Math.\/}, {\bf 89\/} (1987), 37--118.

\bibitem{22} {\sc Harder, G.:}  Some results on the Eisenstein cohomology of
arithmetic subgroups of $GL_n$~.  In: J.P. Labesse, J. Schwermer (Eds):
{\em Cohomology of arithmetic groups and automorphic forms, Lect. Not.
Math.\/}, {\bf 1447\/} (1993), 85--153.

\bibitem{23} {\sc Haagerup, U., Bisch, D.:}  Composition of subfactors: new examples of infinite depth subfactors.  {\em Ann. Scient. \'Ecole Norm. Sup.\/}, {\bf 29\/} (1996), 329--383.

\bibitem{24} {\sc Jones, V.F.R.:} Index for subfactors.  {\em Invent. Math.\/}, {\bf 72\/} (1983), 1--25.

\bibitem{25} {\sc Kasparov, G.G.:}  The operator $K$-functor and extensions
of $C^*$-algebras.  {\em Math. USSR Invest.\/}, {\bf 16\/} (1981),
513--572.

\bibitem{26} {\sc Kasparov, G.G.:}  Equivariant $KK$-theory and the Novikov
conjecture.  {\em Invent. Math.\/}, {\bf 91\/} (1988), 147--201.

\bibitem{27} {\sc Lang, S.:}  Algebraic number theory.  Addison Wesley, 1968.

\bibitem{28} {\sc Murray, F.J., von Neumann, J.:}  On rings of operators. 
{\em Annals of Math.\/}, {\bf 37\/} (1936), 116--129.

\bibitem{29} {\sc Pierre, C.:}   $n$-dimensional global correspondences of Langlands.  Preprint Arxiv.org (2005), RT/0510348.

\bibitem{30} {\sc Pierre, C.:}   Introducing bisemistructures.  Preprint Arxiv.org (2006), GM/0607624.

\bibitem{xx} {\sc Pierre, C.:}   Algebraic quantum theory.  Preprint Arxiv.org (2004), math-ph/0404024.

\bibitem{xxx} {\sc Pierre, C.:}   $n$-dimensional geometric shifted global bilinear correspondences of Langlands on mixed motives - III.  Preprint arxiv.math RT/0709.3383 v3 (2009).

\bibitem{31} {\sc Powers, R.T.:}  Representations of uniformly hyperfinite
algebras and their associated von Neumann rings.  {\em Annals of
Math.\/}, {\bf 86\/} (1967), 138--171.

\bibitem{32} {\sc Schwermer, J.:}  Cohomology of arithmetic groups,
automorphic forms and $L$-functions.  In: Labesse, J.P., Schwermer, J.
(Eds.): {\em Cohomology of arithmetic groups and automorphic forms, Lect.
Notes Math.\/}, {\bf 1447\/} (1990), pp. 1--29.

\bibitem{33} {\sc Serre, J.P.:}  Groupes de Galois sur $Q$~.  {\em Sém.
Bourbaki\/} {\bf 689\/} (1987), 73--85.


\bibitem{34} {\sc Swan, R.:}  Vector bundles and projective modules.  {\em
Trans. Amer. Math. Soc.\/}, {\bf 105\/} (1962), 264--277.

\bibitem{35} {\sc Winter, D.:}  The structure of fields.  {\em Grad. Texts
in Math.\/}, {\bf 16\/} (1974), Springer Verlag.

\end{thebibliography}
